\documentclass[review]{siamonline220329}
\usepackage{siunitx}
\usepackage{amsfonts} 					
\usepackage{amsmath}			
\usepackage{stmaryrd}					
\usepackage{amssymb}					
\usepackage{cleveref}
\usepackage{overpic}
\usepackage{import}

\usepackage{todonotes}

\newcommand{\dphi}{\dot{\phi}}
\newcommand{\ddphi}{\ddot{\phi}}
\newcommand{\dpsi}{\dot{\psi}}
\newcommand{\ddpsi}{\ddot{\psi}}
\newcommand{\dtheta}{\dot{\theta}}
\newcommand{\ddtheta}{\ddot{\theta}}

\newcommand{\dx}{\dot{x}}
\newcommand{\dy}{\dot{y}}
\newcommand{\dz}{\dot{z}}
\newcommand{\ddx}{\ddot{x}}
\newcommand{\ddy}{\ddot{y}}
\newcommand{\ddz}{\ddot{z}}

\newcommand{\dPsi}{\dot{\Psi}}
\newcommand{\dPhi}{\dot{\Phi}}
\newcommand{\dTheta}{\dot{\Theta}}
\newcommand{\order}[1]{\mathcal{O}\left(#1\right)}

\newcommand{\ttwo}{\tau_2}
\DeclareMathOperator\sign{sign}
\DeclareMathOperator\arctanh{arctanh}
\newcommand{\diff}[2][t]{\frac{\mathrm{d} #2}{\mathrm{d}#1}}

\newcommand{\primediff}[1]{#1^\prime}
\newcommand{\primeddiff}[1]{#1^{\prime\prime}}
\newcommand{\Bzero}{B_0}
\newcommand{\Btwo}{B_{-2}}

\renewcommand{\vec}[1]{\underline{#1}}

\newcommand{\mat}[1]{\mathrm{#1}}

\def\app#1#2{%
  \mathrel{%
    \setbox0=\hbox{$#1\sim$}%
    \setbox2=\hbox{%
      \rlap{\hbox{$#1\propto$}}%
      \lower1.1\ht0\box0%
    }%
    \raise0.25\ht2\box2%
  }%
}

\title{The rocking can: a reduced equation of motion and a matched asymptotic solution\thanks{Submitted to the editors 2 Feb 2023.
\funding{SJH would like to thank the Hungarian Academy of Sciences for support through its Distinguished Guest Scientist Programme.}}}

\author{B. W. Collins \and C. L. Hall
\and S. J. Hogan\thanks{Department of Engineering Mathematics, University of Bristol, Bristol BS8 1TW, United Kingdom (\email{ben.collins@bristol.ac.uk, cameron.hall@bristol.ac.uk, s.j.hogan@bristol.ac.uk}). 
Corresponding Author: ORCiD:  0000-0001-6012-6527}}

\headers{Rocking can}{B. W. Collins, C. L. Hall, and S. J. Hogan}

\date{Received: date / Accepted: date}

\begin{document}

\maketitle

\begin{abstract}
The rocking can problem \cite{Srinivasan2008} consists of a empty drinks can standing upright on a horizontal plane which, when tipped back to a single contact point and released, rocks down towards the flat and level state. At the bottom of the motion, the contact point moves quickly around the rim of the can. The can then rises up again, having rotated through some finite \textit{angle of turn} $\Delta\psi$. 
We recast the problem as a second order ODE and find a Frobenius solution. We then use this Frobenius solution to derive a reduced equation of motion. 
The rocking can exhibits two distinct phenomena: behaviour very similar to an inverted pendulum, and dynamics with the angle of turn.
This distinction allows us to use matched asymptotic expansions to derive a uniformly valid solution that is in excellent agreement with numerical calculations of the reduced equation of motion. The solution of the inner problem was used to investigate of the angle of turn phenomenon. 
We also examine the motion of the contact locus $\vec{x}_l$ and see a range of different trajectories, from circular to petaloid motion and even cusp-like behaviour. 
Finally, we obtain an approximate lower bound for the required coefficient of friction to avoid slip. 
\end{abstract}

\begin{keywords} 
Rocking can, Frobenius solution, matched asymptotic expansion
\end{keywords}

\section{Introduction}

\begin{figure}
        \centering{}
        \begin{overpic}[percent,width=0.9\textwidth]{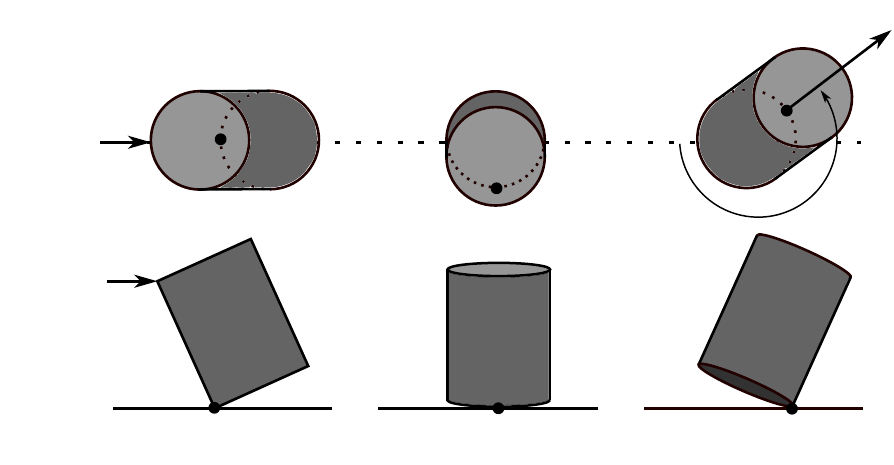}
            \put(15,50){a)}
            \put(48,50){b)}
            \put(80,50){c)}

            \put(0,45){top view}
            \put(0,9){side view}
            \put(93,30){$\Delta\psi$}

            \put(10,33){push}
            \put(10,22){push}
        \end{overpic}
	\caption{The angle of turn phenomenon. a) The can is tilted about a point on the rim given by the black circle. On release, the can falls down. b) As the can approaches the flat state the contact point rapidly races around the rim of the can. c) The can rises up again, pivoting about the contact point. The contact point has moved through an angle $\Delta\psi$ around the rim of the can.}
	\label{fig01:can_fall}
\end{figure}
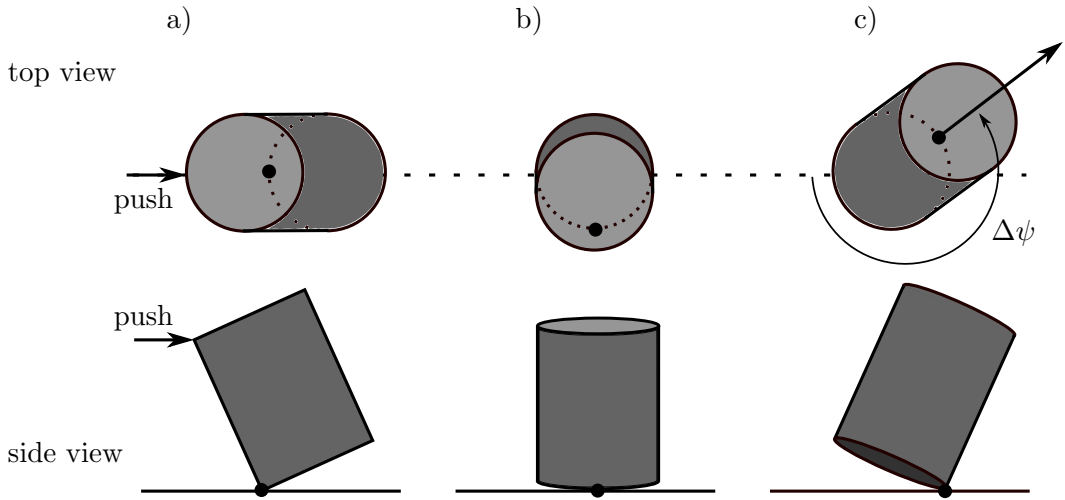

Take an empty drinks can, place it on a horizontal, hard surface, and balance it about a point on its rim (\cref{fig01:can_fall}). When the can is released, with a gentle push along the centre line, it rocks downwards like an inverted compound pendulum, rotating about the contact point. When the can is almost flush with the surface, it appears to bounce and then rocks back up again. During the ``bounce'' phase, the contact point moves rapidly around the rim of the can. When the can rocks back up, the contact point is not diametrically opposed to the starting direction.  

Experiments by Srinivasan and Ruina \cite{Srinivasan2008} determined that the can rotates through an \textit{angle of turn} $\Delta\psi\approx\pm217^\circ$ (\cref{fig01:can_fall}), where the sign is determined by initial conditions. 
Employing small angle approximations and formal assumptions on the dynamics, they estimated $\Delta\psi=\pm202^\circ$. Further analysis \cite{Srinivasan2009} showed that a small off-centre point mass results in chaotic motion of the can.

In earlier work, Cushman and Duistermaat \cite{Cushman2006} had studied the nearly flat falling motions of a thin disk and uncovered similar behaviour. Batista \cite{Batista2008} considered a thick disk, inclined at small angles. Collisions with the surface were shown to be dependent upon two constants of integration, but the physical meanings of these constants were not explored. 

In this paper we study the rocking can problem \cite{Srinivasan2008} from an asymptotic perspective. We consider a rotationally symmetric can rolling and spinning on a horizontal plane with coefficient of friction $\mu$. Particular attention is paid to the motion with small  angular momenta.

In \cref{sec01:derivation} we rederive the equations of motion \cite{Srinivasan2008} and non-dimensionalise them. In \cref{sec01:steady_motions} we determine the static equilibria and steady motions in the problem, along with their stability properties. In \cref{sec01:reduction}, we show that the governing equations can be reduced to a single second order ODE \cref{eq01:frob_ode}. We exploit the presence of a regular singular point in \cref{eq01:frob_ode} to derive a Frobenius solution, which in turn can be used to reduce the equations of motion to one singularly perturbed planar ODE \cref{eq01:Phi_reduced}.
In \cref{sec01:asymptotic_analysis} we carry out an asymptotic analysis of \cref{eq01:Phi_reduced}, gain a uniformly valid approximation for the dynamics and make rigorous the formal assumptions of \cite{Srinivasan2008}. 
\Cref{sec01:physical_phenomena} contains analysis of some properties of the rocking can phenomenon. We calculate the angle of turn and obtain the same expression as \cite{Srinivasan2008}. We also find the condition for the can to fall either clockwise or anticlockwise, extending the work of \cite{Cushman2006}. In addition, we find that the contact locus can move in a circle at variable speed, as well as in cusp-like and petaloid patterns.
 Finally, we test the feasibility of the angle of turn phenomenon by computing a lower bound for the coefficient of friction. 

\section{Derivation of equations of motion}\label{sec01:derivation}

The equations of motion for a can rolling on a rough horizontal plane have been derived \cite{Borisov2017, Leine2009, Ma2014, Srinivasan2008}. In this section, we establish our notation and  rederive equations in the manner of Srinivasan and Ruina \cite{Srinivasan2008}.

\begin{figure}
	\centering
	\begin{overpic}[width =\textwidth,percent]{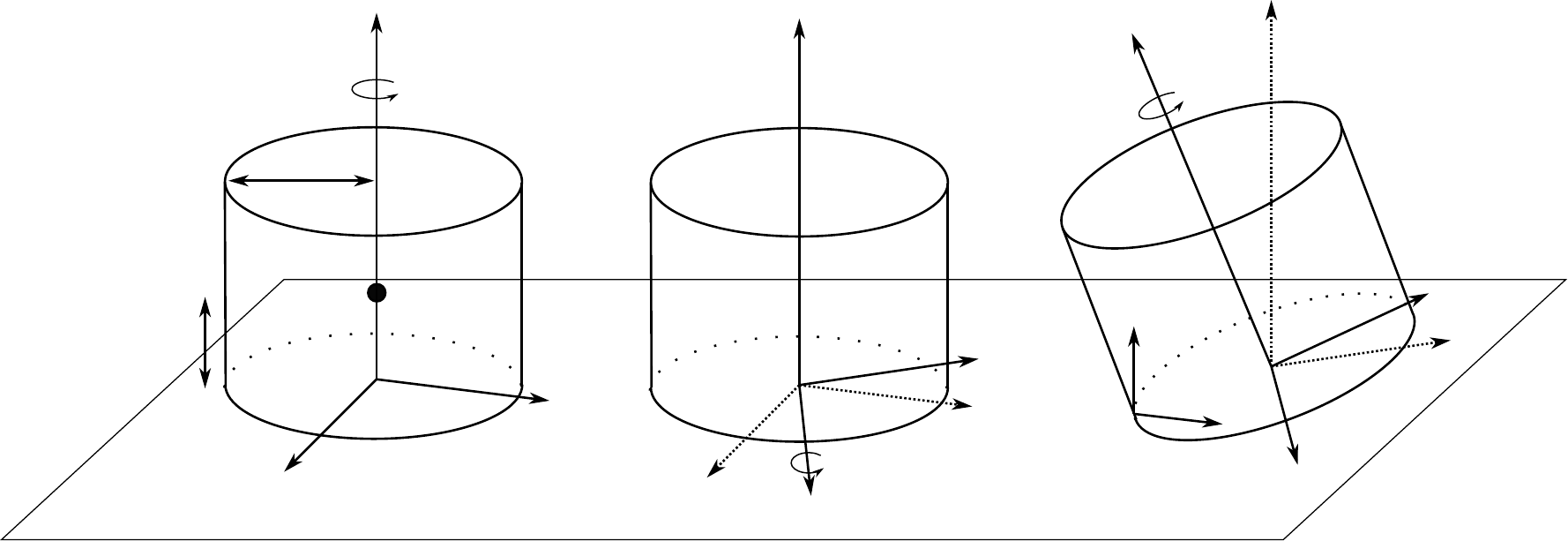}
\put(10,1){$\mathcal{G}$}
\put(40,1){$\mathcal{I}$}
\put(70,1){$\mathcal{B}$}

\put(19,23.5){$R$}
\put(10,12){$H$}
\put(26,14.3){$\vec{G}$}

\put(25,30){$z^\mathcal{G}$}
\put(16,4.5){$x^\mathcal{G}$}
\put(32,6){$y^\mathcal{G}$}
\put(20,28){$\dpsi$}

\put(52,30){$z^\mathcal{I}$}
\put(54,4){$x^\mathcal{I}$}
\put(62,12.5){$y^\mathcal{I}$}
\put(61,9.4){$\psi$}
\put(48.5,2){$\dphi$}

\put(69,28){$z^\mathcal{B}$}
\put(79.25,3.75){$x^\mathcal{B}$}
\put(90.5,17.5){$y^\mathcal{B}$}
\put(72,14.6){$\vec{N}$}
\put(78,9){$\vec{F}$}
\put(71,5){$\vec{P}$}
\put(88,13){$\phi$}
\put(76,29){$\dtheta$}
	\end{overpic}
	\caption{The three reference frames; global $\mathcal{G}$, intermediate $\mathcal{I}$ and body frame  $\mathcal{B}$, given by successive rotations $\psi$ and $\phi$.}
	\label{fig01:can_set_up}
\end{figure}
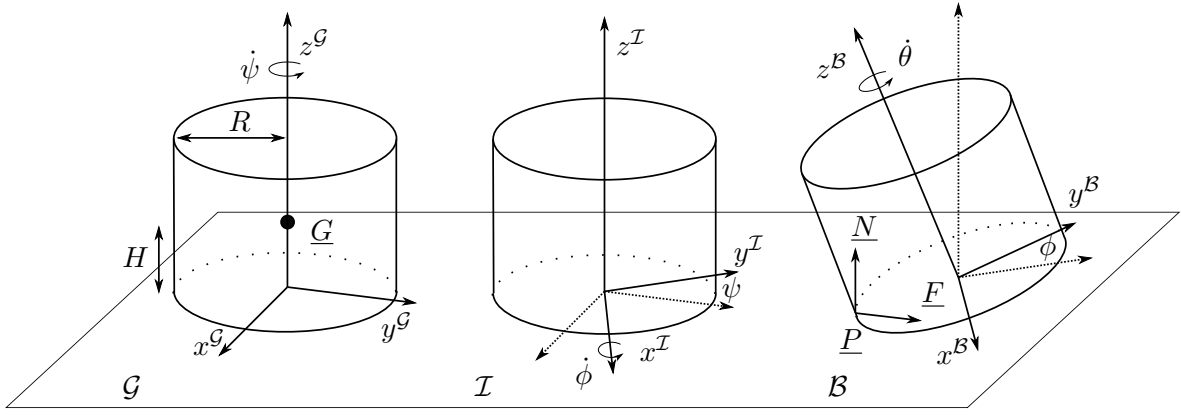

The can, shown in \cref{fig01:can_set_up}, is a rigid, rotationally symmetric cylinder of mass $m$, with height $2H$, radius $R$ and moment of inertia tensor $\mat{I} =\text{diag}\{A,A,C\}$, where $C$ is the moment of inertia about the symmetry axis and $A$ the moment of inertia about the non-symmetry axes. The can moves on a rough horizontal plane with a coefficient of friction $\mu$, that is assumed large enough to ensure rolling motion. A normal reaction $\vec{N}$ and friction force $\vec{F}$ act at the contact point $\vec{P}$.

To describe the orientation of the cylinder, we require three reference frames: the global frame $\mathcal{G}$, an intermediate frame $\mathcal{I}$, and the body frame $\mathcal{B}$, see \cref{fig01:can_set_up}. The frames are defined by Euler angles. In $\mathcal{G}$, axes are aligned with the horizontal plane. Rotation by the \textit{precession} angle, $\psi$, around the $z^\mathcal{G}$ axis gives $\mathcal{I}$. Subsequent rotation by the \textit{nutation} angle, $\phi$, about the $y^\mathcal{I}$ axis brings the can into $\mathcal{B}$.  We also require a final rotation $\theta$, the \textit{rotation} angle, about the $z^\mathcal{B}$ axis. But since this axis is aligned with the symmetry axis of the cylinder, another frame is not required.

The $3\times3$ rotation matrices converting frame $i$ to frame $j$ are given by $\mat{R}_{ij}$, where
\begin{align}
	\mat{R}_{\mathcal{G}\mathcal{I}} = \begin{pmatrix}
		\cos\psi & \sin\psi &0\\
		-\sin\psi & \cos\psi & 0 \\
		0 & 0 & 1 
	\end{pmatrix},\quad 	\mat{R}_{\mathcal{I}\mathcal{B}} = \begin{pmatrix}
		\cos\phi & 0 & \sin\phi \\
		0 & 1 & 0 \\
		-\sin\phi &  0  &\cos\phi
	\end{pmatrix}. \label{eq01:rotation_matrices}
\end{align}
The equations of motion are given by 
\begin{align}
	m\vec{a}_{G}^\mathcal{G} &=  -mg\hat{\vec{z}}^\mathcal{G}  + \vec{N}^\mathcal{G} + \vec{F}^\mathcal{G},\label{eq01:n2l_lin}\\
	\mat{I}\dot{\vec{\Omega}}^\mathcal{B}  + \vec{\omega}^\mathcal{B}\times \mat{I}\vec{\Omega}^\mathcal{B} &=\vec{GP}^\mathcal{B}\times\mat{R}_{\mathcal{G}\mathcal{B}}\left(\vec{N}^\mathcal{G}+ \vec{F}^\mathcal{G}\right),\label{eq01:n2l_rot}
\end{align}
where the superscripts indicate the reference frame of the vector. Note that the force balance (\ref{eq01:n2l_lin}) is expressed in the global frame, whereas the moment balance (\ref{eq01:n2l_rot}) is expressed in the body frame, because the moment of inertia tensor $\mat{I}$ is aligned with the can. $\vec{a}_{G}^\mathcal{G}$ is the acceleration of the centre of mass in $\mathcal{G}$. $\hat{\vec{z}}^\mathcal{G}$ denotes the unit vector parallel to the $z^\mathcal{G}$ axis. $\vec{GP}^\mathcal{B} = (-R,0,-H)^\intercal$ is the vector from the centre of mass $\vec{G}$ to the contact point $\vec{P}$. $\vec{\Omega}^\mathcal{B} = (\dpsi\sin\phi, -\dphi,\dpsi\cos\phi +\dtheta)^\intercal$ is the angular velocity vector in $\mathcal{G}$. $\vec{\omega}^\mathcal{B} = (\dpsi\sin\phi, - \dphi,\dpsi\cos\phi)^\intercal $ is the angular velocity of $\mathcal{B}$ about $\mathcal{G}$ and $\mat{R}_{\mathcal{G}\mathcal{B}}=\mat{R}_{\mathcal{G}\mathcal{I}}\mat{R}_{\mathcal{I}\mathcal{B}}$.

If we assume that the can is rolling without slipping then $\vec{v}_{P}^\mathcal{G}$, the velocity of the contact point $\vec{P}$  in $\mathcal{G}$, is zero and so
\begin{align}
	\vec{v}_{P}^\mathcal{G} = \vec{v}_{G}^\mathcal{G} + \mat{R}_{\mathcal{B}\mathcal{G}}\left(\vec{\Omega}^\mathcal{B}\times \vec{GP}^\mathcal{B}\right) = \vec{0}. \label{eq01:constraint} 
\end{align} 
Hence the velocity of the centre of mass $\vec{v}_{G}^\mathcal{G}$ is given by 
\begin{align}
	\vec{v}_{G}^\mathcal{G}=\begin{pmatrix}
		\dx_G\\ \dy_G \\ \dz_G
	\end{pmatrix}^\mathcal{G} = \begin{pmatrix}
		-\cos\psi(\dphi(R\sin\phi + H\cos\phi) - \sin\psi( \dpsi(R\cos\phi - H\sin\phi) +  R\dtheta)\\
		\cos\psi(\dpsi(R\cos\phi  - H\sin\phi) + R\dtheta) - \sin(\psi)(\dphi(R\sin\phi + H\cos\phi))\\
		\dphi(R\cos\phi - H\sin\phi)
	\end{pmatrix}^\mathcal{G}. \label{eq01:com_velocity}
\end{align}
Differentiating the velocities in \cref{eq01:com_velocity} yields the acceleration of the centre of mass $\vec{a}_{G}^\mathcal{G}=(\ddx_G^\mathcal{G}, \ddy_G^\mathcal{G}, \ddz_G^\mathcal{G})^\intercal$ given by
\begin{subequations}\label{eq01:lin_eom}
\begin{align}
	\ddx_G^\mathcal{G} =& L(\phi,\dphi,\ddphi,\dpsi,\dtheta)\cos\psi+M(\phi,\dphi,\dpsi,\ddpsi,\ddtheta)\sin\psi,\\
	\ddy_G^\mathcal{G} =& M(\phi,\dphi,\dpsi,\ddpsi,\ddtheta)\cos\psi + L(\phi,\dphi,\ddphi,\dpsi,\dtheta)\sin\psi,\\
	\ddz_G^\mathcal{G} =& \ddphi(R\cos\phi - H\sin\phi) - \dphi^2(R\sin\phi + H\cos\phi)
\end{align}
\end{subequations}
where
\begin{subequations}
    \begin{align}
        L(\phi,\dphi,\ddphi,\dpsi,\dtheta)=&-(R\sin\phi + H\cos\phi)\ddphi - (R\cos\phi - H\sin\phi)\dpsi^2
  - R\dtheta\dpsi\\
  &- (R\cos\phi - H\sin\phi)\dphi^2\\
M(\phi,\dphi,\dpsi,\ddpsi,\ddtheta)=&(R\cos\phi - H\sin\phi)\ddpsi + R\ddtheta - 2(R\sin\phi +H\cos\phi)\dphi\dpsi.
    \end{align}
\end{subequations}
Substituting \cref{eq01:lin_eom} into the force balance \cref{eq01:n2l_lin} determines the normal and friction forces
\begin{align}
	F_x =& m\left(L(\phi,\dphi,\ddphi,\dpsi,\dtheta)\cos\psi 
	 +M(\phi,\dphi,\dpsi,\ddpsi,\ddtheta)\sin\psi\right),\label{eq01:fx}\\
	F_y = &m\left(M(\phi,\dphi,\dpsi,\ddpsi,\ddtheta)\cos\psi
 + L(\phi,\dphi,\ddphi,\dpsi,\dtheta)\sin\psi\right),\label{eq01:fy}\\ 
	N =& mg+ m\ddphi(R\cos\phi - H\sin\phi) - m\dphi^2(R\sin\phi + H\cos\phi),\label{eq01:normal}
\end{align}
where $F_x$ and $F_y$ are the components of the friction force $\vec{F}$ along the $x^\mathcal{G}$ and $y^\mathcal{G}$ axes. 

The scalar normal force $N$ is given by $\vec{N} = N\hat{\vec{z}}^\mathcal{G}$. 
Hence \cref{eq01:n2l_rot} become \cite{Srinivasan2008}
\begin{subequations}\label{eq01:raw_eom}
\begin{align}
	   &\left((A + mH^2)\sin\phi - mHR\cos\phi\right)\ddpsi  - mHR\ddtheta = (C -2A-2mH^2)\dpsi\dphi\cos\phi\label{eq01:psi_ugly}\\
       &\quad\quad\quad+ C\dphi\dtheta - 2mHR\dpsi\dphi\sin\phi,\nonumber\\
		&\vspace{1cm}\left(mR^2 + mH^2 + A\right)\ddphi =   \left((A + mH^2 - C- mR^2)\sin\phi\cos\phi - mRH\cos(2\phi)\right)\dpsi^2\label{eq01:phi_ugly}\\
		&\quad\quad\quad - mg(R\cos\phi - H\sin\phi)-((C+mR^2)\sin\phi + mRH\cos\phi)\dtheta\dpsi,\nonumber\\
	   &\vspace{1cm}\left((C + mR^2)\cos\phi - mRH\sin\phi)\ddpsi + (C+mR^2\right)\ddtheta =C\dpsi\dphi\sin\phi + 2mR\dpsi\dphi(R\sin\phi  + H\cos\phi).\label{eq01:theta_ugly}
\end{align}
\end{subequations}
\Cref{eq01:raw_eom,eq01:lin_eom} form the equations of motion for the can, with \cref{eq01:lin_eom} being cyclic. We rescale the lengths\footnote{The choice of $R$, rather than $H$, to non-dimensionalise lengths avoids large quantities when considering thin disks where $0<H\ll R$. } by $R$ and the moments of inertia by $mR^2$, introducing

\begin{align}
		h = \frac{H}{R},\quad x =\frac{x_G}{R}, \quad y = \frac{y_G}{R},\quad z = \frac{z_G}{R},\quad a = \frac{A}{mR^2},\quad c = \frac{C}{mR^2}.\label{eq01:scalings}
\end{align}
We scale $t$ by $\sqrt{R/g}$ and overload the notation so that the dot notation means differentiation with respect to the scaled time.

\begin{figure}
	\centering
	\begin{overpic}[percent,width=0.9\textwidth]{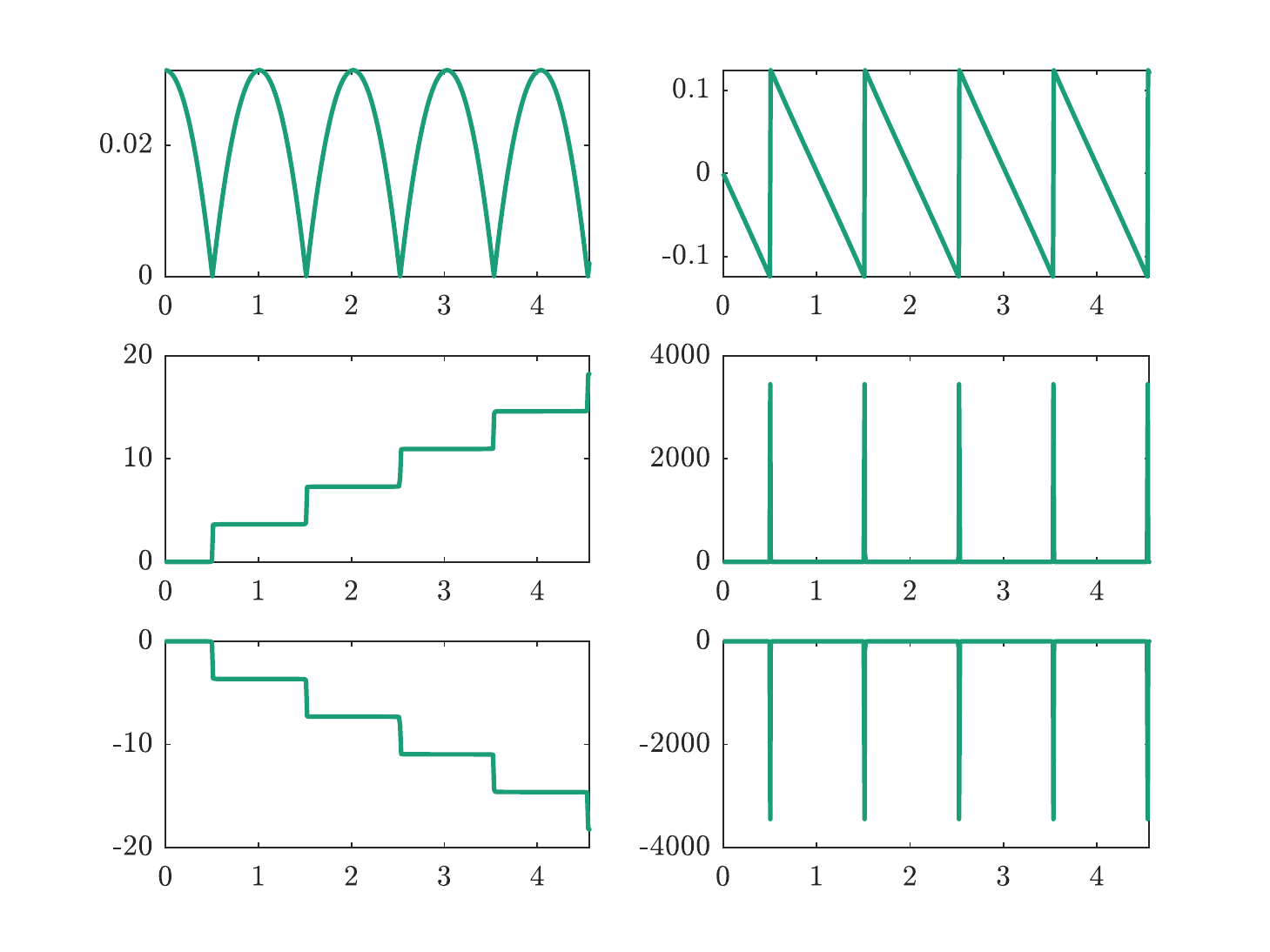}
			\put(5,62){$\phi$}
			\put(5,38){$\psi$}
			\put(5,16){$\theta$}
			\put(48,62){$\Phi$}
			\put(48,38){$\Psi$}
			\put(48,16){$\Theta$}
			\put(28,48){$t$}
			\put(28,26){$t$}
			\put(28,3){$t$}
			\put(72,48){$t$}
			\put(72,26){$t$}
			\put(72,3){$t$}
	\end{overpic}
	\caption{Numerical solution of the non-dimensionalised equations of motion \cref{eq01:non_dimensionalised_eom}. The material parameters and initial conditions, used throughout this paper, are given in \cref{eq01:material_parameters,eq01:initial_conditions} respectively.
 }\label{fig01:non_dim_solution}
\end{figure}
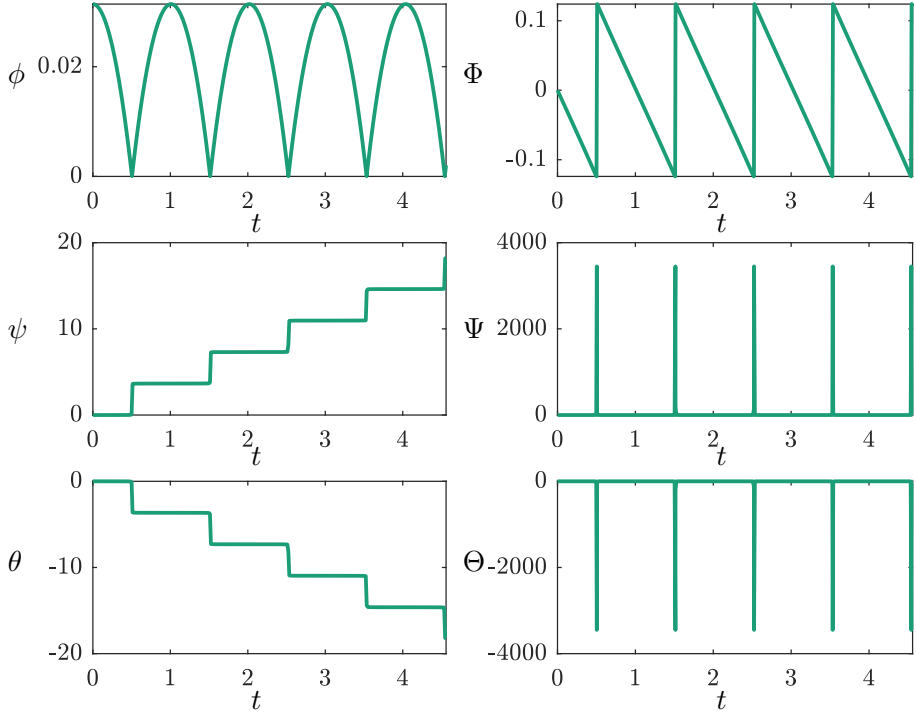

\Cref{eq01:raw_eom} then becomes
\begin{subequations}\label{eq01:non_dimensionalised_eom}
\begin{align}
		\dPsi\sin\phi  =& kc_p\Phi\Theta +  ((kc_p-2)\cos\phi - hk\sin\phi ) \Phi\Psi,\label{eq01:Psi}\\
		\dTheta\sin\phi =& (-kc_p\cos\phi + hk\sin\phi)\Phi\Theta  + ( -kc_p\cos^2\phi - ka_p\sin^2\phi + hk\sin 2\phi + 2)\Phi\Psi,\label{eq01:Theta}\\
		\dPhi(a_p+ 1) =& ((a_p - c_p)\sin\phi\cos\phi - h\cos2\phi)\Psi^2 - (c_p\sin\phi + h\cos\phi)\Theta\Psi + ( h\sin\phi - \cos\phi),\label{eq01:Phi}\\ 
  \dpsi=&\Psi,\label{eq01:psi}\\
  \dtheta=&\Theta,\label{eq01:theta}\\
		\dphi =& \Phi, \label{eq01:phi}
\end{align}
\end{subequations}
where the constants
\begin{align}
 	a_p = a + h^2,\quad c_p = c + 1
 \end{align}  
 are the scaled moments of inertia $a$ and $c$ about the contact point $\vec{P}$, and we set
 \begin{align}
 	k = \frac{c}{a + ca_p}.
 \end{align}
  Throughout this paper, we consider the can \cite{Srinivasan2008} to have mass $m = \SI{4.3e-2}{\kilogram}$, height $H =\SI{5.45e-2}{\meter}$ and radius $R = \SI{3.7e-2}{\meter}$. The moments of inertia\footnote{Our parameters differ to those used by Srinivasan and Ruina \cite{Srinivasan2008}} are $A = \SI{6.97e-5}{\kilogram\meter\squared}$ and $C = \SI{5.89e-5}{\kilogram\meter\squared}$, and 
  \begin{align}
      a = 0.727, \quad c = 0.615, \quad a_p = 2.897, \quad c_p = 1.6156, \quad  h = 1.473, \quad  k = 0.245. \label{eq01:material_parameters}
  \end{align}
  Numerical solutions of \cref{eq01:non_dimensionalised_eom}, performed in \textsc{matlab} using \textsc{ode15s} to cope with the stiff ODEs, are shown in \cref{fig01:non_dim_solution}; compare \cite[Fig. 5]{Srinivasan2008}. Fast changes in the variables can be observed when $\phi$ is small. The angle of turn $\Delta\psi$ is the step-like change in $\psi$. Unless otherwise stated, throughout the paper we take initial conditions 
  \begin{align}
   (\psi,\dpsi,\phi,\dphi,\theta,\dtheta) = (0,0.1001,\pi/100,0,0,-0.1000). \label{eq01:initial_conditions}   
  \end{align}

The rocking can phenomenon cannot be obtained by setting $\Psi=\Theta=0$, as we now demonstrate. Let us  assume $\Phi\ne0$. \Cref{eq01:non_dimensionalised_eom} then yields an integrable subsystem, where the can is pivoting about its rim, but not spinning or rolling. We find 
\begin{align}
	\dPhi=\ddphi = \frac{h\sin\phi - \cos\phi}{a_p+1}, \label{eq01:nonlin_pendulum}
\end{align}
with solution 
\cite{ZaitsevValentinFandPolyanin2002}
\begin{align}
	\phi(t) = \arctan(1/h) \pm \int^{\tau}(-2\cos s + C_1)^{-1/2}\mathrm{d}s + C_2,
\end{align}
where $\tau = \sqrt{\frac{h^2+1}{a_p +1}}t$ and $C_1, C_2$ are integration constants. This solution passes through $\phi=0$, impacting the plane at all points on its rim \textit{simultaneously}, something that is not observed in experiments \cite{Srinivasan2008}.

\section{Equilibria and steady motions}\label{sec01:steady_motions}

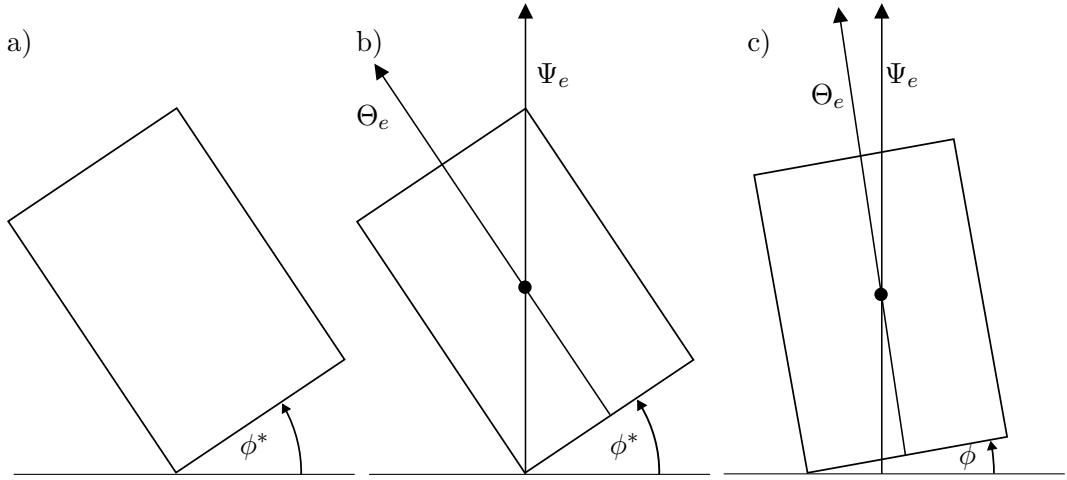
\begin{figure}
		\centering
		\begin{overpic}[width = 0.9\textwidth]{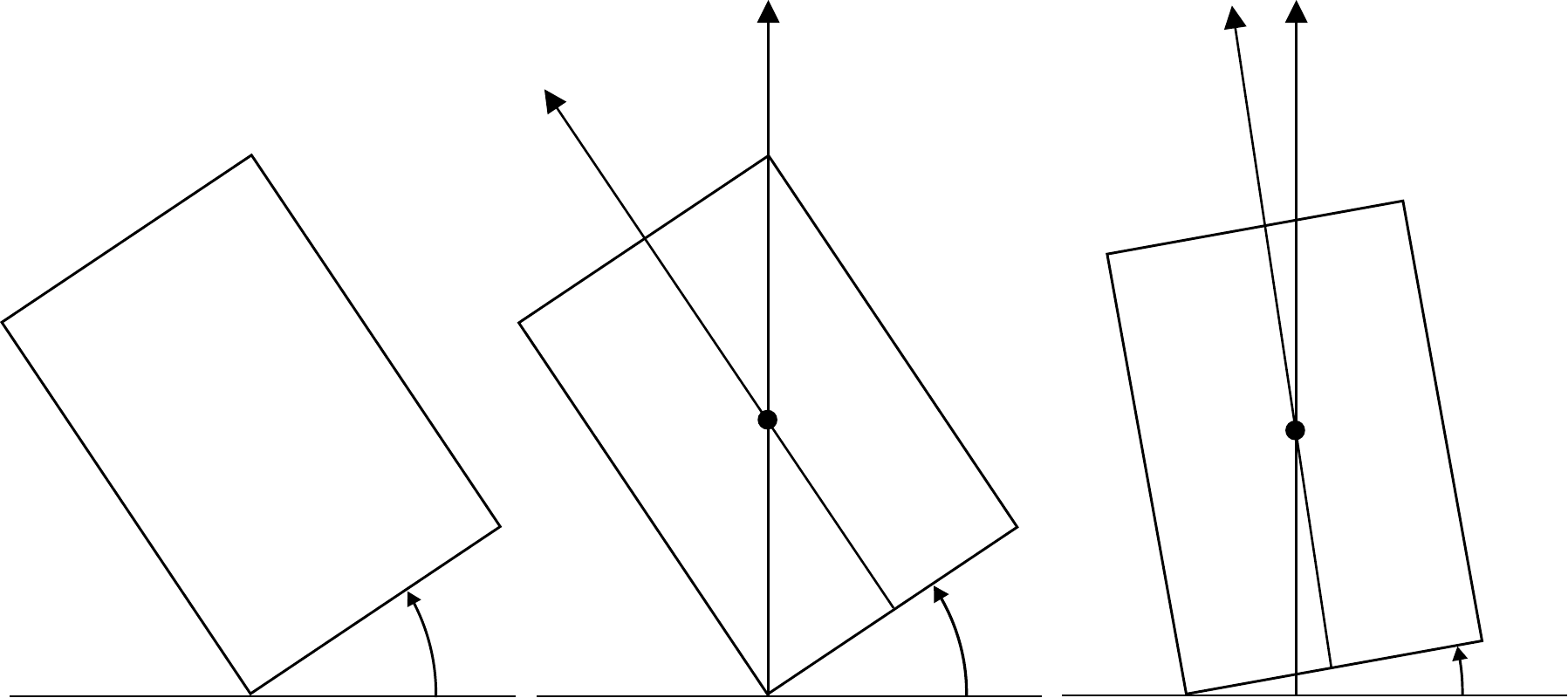}
		\put(0,40){a)}
		\put(22,2){$\phi^*$}

		\put(33,40){b)}
		\put(57,2){$\phi^*$}
		\put(33,33){$\Theta_e$}
		\put(50,37){$\Psi_e$}

		\put(70,40){c)}
		\put(90,1){$\phi$}
		\put(76,35){$\Theta_e$}
		\put(83,37){$\Psi_e$}

		\end{overpic}
	\caption{a) $S_\text{static}$ \cref{eq01:static_eqm}: the can is static and balanced on its rim at $\phi=\phi^*$ \cref{eq:phistar}. b)  $S_\text{bal}$ \cref{eq01:bal_eqm}: the can is balanced on its rim, $\phi=\phi^*$, rolling with constant $\Theta_e, \Psi_e$. c) $S_\text{steady}$ \cref{eq01:equilibrium}: the can is rolling with constant $\Theta_e, \Psi_e$, $\phi\ne\phi^*$}. \label{fig01:eqm_comparison}
\end{figure}

In this section we locate the equilibria and steady motions of  \cref{eq01:non_dimensionalised_eom} and determine their stability. Physically, we expect an equilibrium $S_\text{static}$ when the can is balanced on its rim and stationary
\begin{align}
	S_\text{static} := \left\{(\Psi,\Theta,\Phi,\phi)= (0,0,0,\phi^*)\right\},\label{eq01:static_eqm}
\end{align}
where the \textit{balancing angle}\footnote{For our can, $\phi^* \approx 34^\circ$. }
\begin{equation}\label{eq:phistar}
    \phi^* := \arctan\left(1/h\right).
\end{equation}
Here all the velocities are zero and the can is tilted at $\phi=\phi^*$, as in \cref{fig01:eqm_comparison}a. $S_\text{static}$ shares similarities with a planar inverted compound pendulum. It is the only static equilibrium of the system.

Now assume constant $\Psi =\Psi_e, \Theta =\Theta_e$, and keep $\Phi = 0,\phi = \phi^*$ as before. Then the set of all steady motions with the can balanced at $\phi = \phi^*$, as in \cref{fig01:eqm_comparison}b, is $S_\text{bal}$ where
\begin{align}
	S_\text{static} \subset S_\text{bal} := \left\{(\Psi,\Theta,\Phi,\phi)= (\Psi_e,\Theta_e,0,\phi^*)\right\}, \label{eq01:bal_eqm}
\end{align}
and $\Psi_e \sin\phi^*\left((a-c)\cos\phi^*\Psi_e - c\Theta_e\right) = 0$ from \cref{eq01:Phi}. 

If $\Psi_e = 0, \Theta_e \ne 0$, then $S_\text{bal}$ represents the balanced can rolling steadily in a straight line, with two non-zero eigenvalues
\begin{align}
	\lambda_{\pm} = \pm \sqrt{\frac{\sqrt{h^2 + 1} - \Theta_e^2c_phk(h + c_p)}{a_p + 1}},\label{eq01:eigenvalues}
\end{align}
and two zero eigenvalues.
For
\begin{align}
	\Theta_e^2 > \Theta^2_{crit} := \frac{\sqrt{h^2 + 1}}{hkc_p(h + c_p)},	
\end{align} 
$S_\text{bal}$ is centre-like, otherwise it is saddle-like. Therefore, for $|\Theta_e|>\Theta_{crit}$, rolling motion will persist (an unbounded solution). For a thin disk, this \textit{critical rolling} has been previously been explored \cite{OReilly1996, Paris2002, Przybylska2016}. 

If $\Psi_e\ne0, \Theta_e \ne 0$, then  $S_\text{bal}$ gives a 1-D manifold of spinning and rolling solutions with $\Theta_e = (a-c)\cos\phi^*\Psi_e/c$. 

Steady rolling can also occur for constant $\Psi_e \ne 0, \Theta_e \ne 0$, with $\Phi=0$, $\phi =\phi_0 \ne\phi^*$, \cref{fig01:eqm_comparison}c. 
Then the set of all equilibria of \cref{eq01:non_dimensionalised_eom} is given by 
\begin{align}
	S_\text{bal} \subset S_\text{steady} = \Big\{(\Psi,\Theta,\Phi,\phi) = (\Psi_e,\Theta_e,0,\phi_0) \Big\}, \label{eq01:equilibrium}
\end{align}
where
\begin{align}
	((a_p - c_p)\sin\phi_0\cos\phi_0 - h\cos2\phi_0)\Psi_e^2 - (c_p\sin\phi_0 + h\cos\phi_0)\Theta_e\Psi_e + ( h\sin\phi_0 - \cos\phi_0) =0.
\end{align}

In addition to steady motions at small $\phi$, there exist equilibria for $\phi_0 > \phi^*$ where the can rolls with the symmetry axis almost horizontal\footnote{These equilibria are central to the `mysterious spinning cylinder' explored by Jackson et al. \cite{Jackson2019}, which consists of a slender cylinder undergoing a steady motion close to $\phi \approx \pi/2$. Upon spinning, symbols drawn on the cylinder disappear or appear to remain stationary depending on their position.}. An example of such steady motion is given by $\phi_0 = \pi/2$, $\Psi_e = \Theta_e = \pm \sqrt{h/(c_p-h)}$.

We now show that the can rolls in a circle for all steady motions $S_\text{steady}$ \cite{LeSaux2005, Ma2016a, OReilly1996, Paris2002}. Applying scalings \cref{eq01:scalings}, we obtain from \cref{eq01:com_velocity}: 
 \begin{subequations}\label{eq01:non_dimensionalised_com}
\begin{align}
	\dx =& -(\Phi(\sin\phi_0 + h\cos\phi_0)\cos\psi - (\Psi_e(\cos\phi_0 - h\sin\phi_0) + \Theta_e)\sin\psi,\\
	\dy =& (\Psi_e(\cos\phi_0 - h\sin\phi_0) + \Theta_e)\cos\psi - (\Phi(\sin\phi_0 + h\cos\phi_0))\sin\psi.
\end{align}
\end{subequations}
On $S_\text{steady}$, $\Phi=0$ and so from \cref{eq01:non_dimensionalised_com} 
\begin{subequations}\label{eq01:comcoords}
\begin{align}
	x(t) &= \frac{\Psi_e(\cos\phi_0 - h\sin\phi_0) + \Theta_e}{\Psi_e}\cos(\Psi_e t),\\
	y(t) &= \frac{\Psi_e(\cos\phi_0 - h\sin\phi_0) + \Theta_e}{\Psi_e}\sin(\Psi_e t),
\end{align}
\end{subequations}
corresponding to circular motion with radius
\begin{align}
r_\text{circ} = \left | \frac{\Psi_e(\cos\phi_0 - h\sin\phi_0) + \Theta_e}{\Psi_e} \right |.
\end{align} 
If 
\begin{equation}\label{eq01:comrest}
\Psi_e(\cos\phi_0 - h\sin\phi_0) + \Theta_e =0,    
\end{equation} 
then $x=y=0$ \cref{eq01:comcoords} and the centre of mass is at rest. 

 On $S_\text{bal}$, where $\phi_0=\phi^*$, $r_\text{circ}=\left |\frac{\Theta_e}{\Psi_e} \right |$. In this case, if $\Theta_e=0, \Psi_e \ne 0$, the centre of mass is at rest, with $r_\text{circ}=0$. But if $\Psi_e=0, \Theta_e\ne0$ the can rolls in a straight line. 

Steady motions with the centre of mass at rest are called \textit{stationary motions} \cite{LeSaux2005,Mcdonald2002}. A can undergoing stationary motion experiences no friction force, just as the frictionless $\mu=0$  case. Including \cref{eq01:comrest} in \cref{eq01:equilibrium} gives the 1-D set of frictionless orbits 
\begin{align}
	S_\text{rest} := \left\{(\Psi,\Theta,\Phi,\phi) =(\Psi_e,\Theta_e,0,\phi_0 | \Psi_e(\cos\phi_0 - h\sin\phi_0) + \Theta_e =0)\right\} \subset S_{steady},
\end{align}
where
\begin{align}
	\Psi_e^2 = \frac{\cos\phi_0 - h\sin\phi_0}{a\sin\phi_0\cos\phi_0 + ch\sin^2\phi_0},\quad\Theta_e = -\Psi_e(\cos\phi_0 - h\sin\phi_0).
\end{align}
The set $S_\text{rest}$ only exists for $\phi_0 < \phi^*$. 

In the next \cref{sec01:reduction} we consider the dynamics of the (non-cyclic) equations \cref{eq01:non_dimensionalised_eom}. 

\section{Frobenius solution and a  reduced equation of motion}\label{sec01:reduction}

There is structure in the equations of motion \eqref{eq01:non_dimensionalised_eom} that has been exploited to give an exact solution \cite{Batista06}. Our derivation is slightly different. Divide \eqref{eq01:Psi} and \eqref{eq01:Theta} by $\Phi$ to get
\begin{subequations}\label{eq01:closed_subsystem}
\begin{align}
	\diff[\phi]{\Psi}\sin\phi   &= kc_p\Theta +  \left((kc_p-2)\cos\phi - hk\sin\phi \right)\Psi, \label{eq01:psi_diff_phi}\\
	\diff[\phi]{\Theta}\sin\phi &= (-kc_p\cos\phi + hk\sin\phi)\Theta  + (-kc_p\cos^2\phi - ka_p\sin^2\phi + hk\sin 2\phi + 2)\Psi, \label{eq01:theta_diff_phi}
\end{align}
\end{subequations}

\begin{figure}
		\centering
		\begin{overpic}[percent,width=0.7\textwidth]{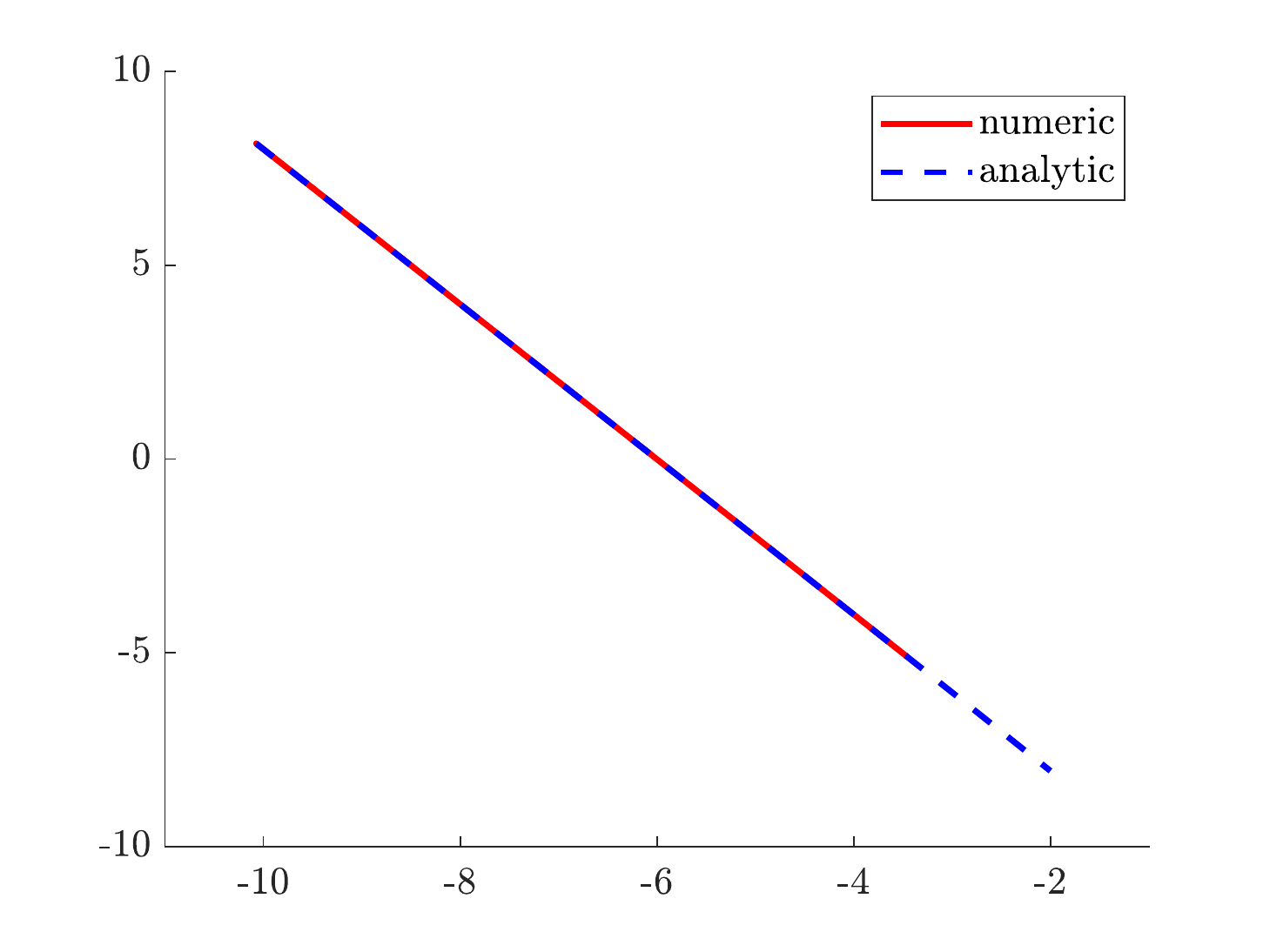}
				\put(-2,40){$\log(\Psi)$}
				\put(50,0){$\log(\phi)$}
		\end{overpic}
	\caption{Comparison of the analytic solution \cref{eq01:hyper_geom} with the numerical solution of the equations of motion \cref{eq01:non_dimensionalised_eom}. The solutions agree and show a clear $\Psi\propto1/\phi^2$ dependency at small $\phi$. 
 The $\Psi^\prime$ initial condition for \cref{eq01:frob_ode} is computed from \cref{eq01:psi_diff_phi}. 
 }\label{fig01:frob_comp}
\end{figure}

Now divide both sides of \cref{eq01:psi_diff_phi} by $\sin\phi$ and differentiate  with respect to $\phi$. Then eliminate $\Theta$ using \cref{eq01:psi_diff_phi} and $\diff[\phi]{\Theta}$ using \cref{eq01:theta_diff_phi} to find
\begin{align}
	\Psi^{\prime\prime} + 3\cot\phi \Psi^\prime - (\gamma  + \beta\cot\phi)\Psi =0, \label{eq01:frob_ode}
\end{align}
 where\footnote{For our can, $\gamma=2.245, \, \beta = 0.361$ from \cref{eq01:material_parameters}.} $\gamma =2+k, \beta = kh$ and $\prime$ denotes differentiation with respect to $\phi$. \Cref{eq01:frob_ode} is a regular, singular ODE with exact solution \cite[eq. (21)]{Batista06}
\begin{align}
 	\Psi(\phi) = & \frac{C_3}{\sin(\phi)^{3/2}}(\cot(\phi) + i)^{- \sqrt{\bar\rho}}(\cot(\phi)-i)^{ \sqrt{\rho}} \times \label{eq01:hyper_geom} \\
  &{_2F_1}\left(\sqrt{\rho} - \frac{1}{2} -\sqrt{\bar\rho}, \sqrt{\rho} + \frac{3}{2} - \sqrt{\bar{\rho}}; 1 - 2 \sqrt{\bar\rho}, \frac{1}{2}- \frac{i}{2} \cot\phi\right)\nonumber  \\
 	&+ \frac{C_4}{\sin(\phi)^{3/2}}(\cot(\phi) + i)^{\sqrt{\bar\rho}}(\cot(\phi)-i)^{\sqrt{\rho}} \times\nonumber  \\  
  &{_2F_1}\left(\sqrt{\rho} + \frac{3}{2} + \sqrt{\bar\rho}, \sqrt{\rho} -\frac{1}{2} + \sqrt{\bar{\rho}}; 1 + 2 \sqrt{\bar\rho}, \frac{1}{2} - \frac{i}{2} \cot\phi\right). \nonumber
 \end{align}
where $\rho = \frac{1}{16}\left(9 + 4\gamma + 4i\beta\right)$, $\bar\rho$ is the complex conjugate, ${_2F_1}$ is the hypergeometric function and $C_3, C_4$ are integration constants. We can find $\Theta(\phi)$ from \cref{eq01:psi_diff_phi,eq01:hyper_geom}.
The exact solution \cref{eq01:hyper_geom} and numerical solution of the full system \cref{eq01:non_dimensionalised_eom} are shown in \cref{fig01:frob_comp} and suggest a strong dependence on $1/\phi^2$ at small $\phi$. Numeric evaluation of the hypergeometric function in \cref{eq01:hyper_geom} with complex parameters and arguments is non-trivial \cite{Doornik2015}. We evaluate \cref{eq01:hyper_geom} with the \textsc{Matlab} function hypergeom().

\subsection{Frobenius series solution}
\Cref{eq01:hyper_geom} is cumbersome to work with. Since we are most interested in the behaviour of the can for small $\phi$, we can solve \cref{eq01:frob_ode} using a Frobenius series \cite{bender78}, since the singular point $\phi=0$ is regular; both $\phi\cot\phi$ and $\phi^2(\gamma + \beta\cot\phi)$ have valid Taylor series expansions. Hence we assume
\begin{align}
	\Psi(\phi) = \sum^\infty_{n=0}a_n\phi^{n+r},
\end{align}
with $a_0 \ne 0$, valid near $\phi=0$. Differentiating the series and substituting into \cref{eq01:frob_ode} gives 
\begin{align}
	\sum a_n(n+r)(n+r-1)\phi^{n+r-2} + 3\cot\phi\sum a_n(n+r)\phi^{n+r-1} - (\gamma + \beta\cot\phi)\sum a_n\phi^{n+r} =0. \label{eq01:frobenius_subs}
\end{align}
Using the series expansion of $\cot\phi$, 
the indicial equation of \cref{eq01:frobenius_subs} is given by
\begin{align}
	r(r+2) = 0.
\end{align}
Therefore $r = r_1 = 0$ and $r=r_2=-2$ and two series solutions $\Psi_0, \Psi_{-2}$ exist for \cref{eq01:frob_ode}. When $r=r_1=0$, the solution takes the form
\begin{align}\label{eq01:Psi0}
	\Psi_0(\phi) = \sum^\infty_{n=0} p_n\phi^n,
\end{align} 
where $p_n$ are constants. Using \cref{eq01:frobenius_subs}, we find 
\begin{align}
	\Psi_0(\phi) = p_0\left(1 + \frac{\beta}{3}\phi + \frac{\beta^2 + 3\gamma}{24}\phi^2 +   \frac{\beta(\beta^2 + 11\gamma)}{360}\phi^3\right) + \order{\phi^4}.
\end{align}
Since $r_1$ and $r_2$ differ by an integer, there will be $\log$ terms \cite{bender78} in the $r=r_2=-2$ solution, which then takes the form
\begin{align}
	\Psi_{-2}(\phi) = C\Psi_0(\phi)\log(\phi) + \sum^\infty_{n=0}q_n\phi^{n-2}.\label{eq0:Psi2}
\end{align}
where $C, q_n$ are constants. Inserting \cref{eq0:Psi2} into \cref{eq01:frob_ode} and simplifying gives
\begin{align}
	&\left[\frac{2C\primediff{\Psi_0}}{\phi} - \frac{C\Psi_0}{\phi^2} + \sum q_n(n-2)(n-3)\phi^{n-4}\right]\\
	&+ 3\cot\phi\left[\frac{C\Psi_0}{\phi} + \sum q_n(n-2)\phi^{n-3}\right]
	- (\gamma + \beta\cot\phi)\left[\sum q_n \phi^{n-2}\right] =0.\nonumber
\end{align}
Using the series expansion of $\cot\phi$, equating orders 
and setting the arbitrary constant $q_2 =0$, we find
\begin{align}
	\Psi_{-2}(\phi) = \frac{-q_0(\beta^2 - \gamma + 2)}{2p_0}\Psi_0(\phi)\log\phi + q_0\left(\frac{1}{\phi^2} - \frac{\beta}{\phi} - \frac{1}{9}\beta\left(-2\beta^2 + 5\gamma - 6\right)\phi\right) + \order{\phi^2}.
\end{align}
Hence the solution to \cref{eq01:frob_ode} for small $\phi$, using the original parameters $\beta=kh$ and $\gamma = 2+k$, is given by
\begin{align}
	\Psi(\phi) =& \Bzero\left(1 + \frac{kh}{3}\phi\right) - \frac{\Btwo((kh)^2 - k)}{2}\left(1 + \frac{kh}{3}\phi + \frac{(kh)^2 + 6 + 3k}{24}\phi^2\right)\log\phi\label{eq01:psi_orig}\\
	 &+ \Btwo\left(\frac{1}{\phi^{2}} - \frac{kh}{\phi} - \frac{1}{9}kh\left(-2(kh)^2 + 4 + 5k\right)\phi\right) + \order{\phi^2},\nonumber
\end{align}
where we have now replaced $p_0$ and $q_0$ by $\Bzero$ and $\Btwo$.
We observe that $\Psi$ has a $\phi^{-2}$ dependence at leading order, as shown numerically in \cref{fig01:frob_comp}.

Recall that $\Theta(\phi)$ and $\Psi(\phi)$ are related through \cref{eq01:psi_diff_phi}. By substituting the series solution for $\Psi$ \cref{eq01:psi_orig} and its derivative, expanding the trigonometric functions and rearranging, we find from \cref{eq01:psi_diff_phi} that
\begin{align}
	\Theta(\phi) =& -\frac{\Btwo}{\phi^2} + \frac{\Btwo kh}{\phi} + \frac{-9\Btwo h^2k^2 + ((4c_p +3)\Btwo - 6\Bzero c_p)k + 12\Bzero - 4\Btwo}{6kc_p}\label{eq01:theta_orig}\\ &+ \Btwo\frac{h^2k-1}{2c_p}\left((c_pk-2) + \frac{kh(c_pk-6)}{3}\phi\right)\log\phi + \order{\phi},\nonumber
\end{align}
Thus to leading order, the series solutions\footnote{The subscript notation $l0,l1$ is used to denote coefficients of $\log$ terms.} for $\Psi$ and $\Theta$ are
\begin{subequations} \label{eq01:series_sols}
\begin{align}
	\Theta(\phi) &= \frac{-\Btwo}{\phi^2} + \frac{\Btwo kh}{\phi} + \Theta_{00} + (\Theta_{l0} + \Theta_{l1}\phi)\log\phi\label{eq01:theta_series},\\
	\Psi(\phi) &= \frac{\Btwo }{\phi^2} + \frac{-\Btwo kh}{\phi} + \Psi_{00} + (\Psi_{l0} + \Psi_{l1}\phi)\log\phi\label{eq01:psi_series}
\end{align}
\end{subequations}
The two expansions suggest that $\Theta + \Psi  = \order{1}$. 
The coefficients $\Theta_{00}, \Psi_{00}, \Theta_{l0}, \Theta_{l1}, \Psi_{l0}, \Psi_{l1}$, which involve $\Bzero, \Btwo$, are given in \cref{app:coefficients_Psi_Theta}.
\Cref{eq01:series_sols} represents an improvement on the relations found by Srinivasan and Ruina \cite[eqs. (7) and (9)]{Srinivasan2008}.

\subsection{Coefficients \texorpdfstring{$\Bzero$}{B0} and \texorpdfstring{$\Btwo$}{B-2}}

\begin{figure}
	\centering
		\begin{overpic}[percent,width=0.7\textwidth]{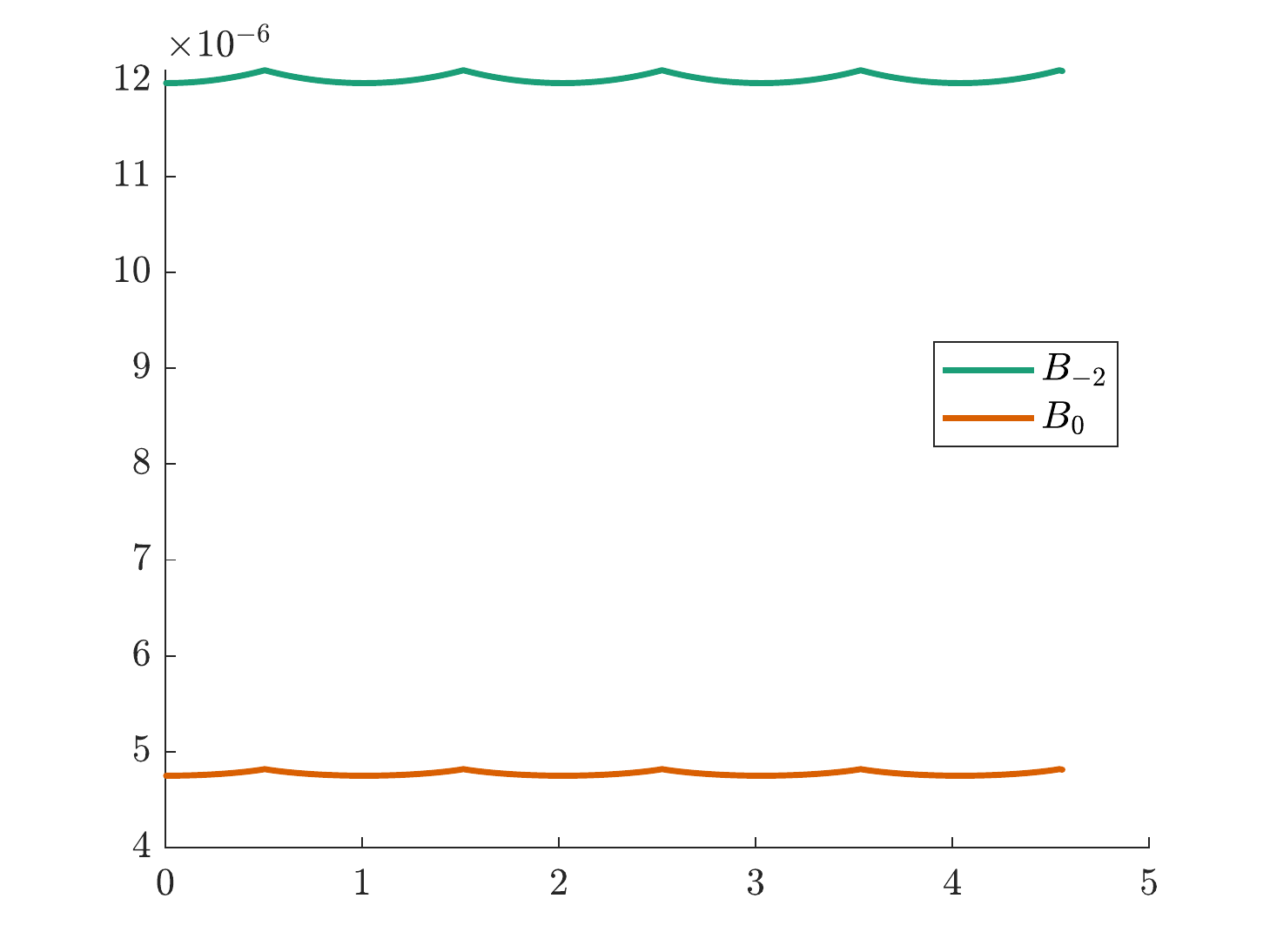}
				\put(3,26){\rotatebox{90}{Conserved quantities}}
				\put(52,2){$t$}
		\end{overpic}
	\caption{The conserved quantities $\Bzero$ and $\Btwo$ computed from the non-dimensionalised equations of motion \cref{eq01:non_dimensionalised_eom}. Small oscillations are visible, due to the truncation of \cref{eq01:A_sol,eq01:B_sol}.
 }\label{fig01:conserved}
\end{figure}

To determine $\Bzero$ and $\Btwo$, we need conditions for \cref{eq01:frob_ode}. But we are given $\Psi_{t=0}$, $\Theta_{t=0}$, $\phi_{t=0}$ and $\Phi_{t=0}$. To determine $\primediff{\Psi}_{t=0}$ we expand \cref{eq01:psi_diff_phi} in small $\phi$ and evaluate at $t=0$
\begin{align}
    \primediff{\Psi}_{t=0} = \frac{kc_p(\Psi_{t=0} + \Theta_{t=0}) - 2\Psi_{t=0}}{\phi_{t=0}}  - hk\Psi_{t=0} + \frac{kc_p(\Theta_{t=0} - 2\Psi_{t=0}) + 4\Psi_{t=0}}{6}\phi_{t=0}+ \order{\phi^3}. \label{eq01:psi_diff_phi_small_phi}
\end{align}
We take \cref{eq01:psi_orig} for $\Psi$, differentiate it with respect to $\phi$ to get a series solution for $\primediff{\Psi}$, evaluate both at $t=0$. We obtain two simultaneous equations in $\Bzero$ and $\Btwo$. 
\begin{subequations}\label{eq01:simultaneous}
\begin{align}
	\Psi_{t=0} &= \Bzero\Psi_0(\phi_{t=0}) + \Btwo\Psi_2(\phi_{t=0}),\\
	\primediff{\Psi}_{t=0} &= \Bzero\primediff{\Psi}_0(\phi_{t=0}) + \Btwo\primediff{\Psi}_2(\phi_{t=0}).
\end{align}
\end{subequations}
We solve for $\Bzero, \Btwo$ to find the series solution
\begin{align}
	\Bzero =& \frac{kc_p}{2}(\Theta_{t=0} + \Psi_{t=0}) - \frac{hk^2c_p}{2}(\Theta_{t=0} + \Psi_{t=0})\phi_{t=0}\label{eq01:A_sol}\\
	& - \frac{h^2k^2-k}{4}(kc_p(\Theta_{t=0} + \Psi_{t=0}) -2\Psi_{t=0})\phi_{t=0}^2 \log(\phi_{t=0})+\nonumber\\
		&+  \frac{1}{12}\left(-kc_p(5\Psi_{t=0} + 2\Theta_{t=0}) + (4 - 3k + 9h^2k^2)\Psi_{t=0}\right)\phi_{t=0}^2 + \order{\phi^3_{t=0}\log(\phi_{t=0})},\nonumber\\
	\Btwo =& \left( - \frac{kc_p}{2}(\Theta_{t=0} + \Psi_{t=0}) + \Psi_{t=0}\right)\phi_{t=0}^2 + \order{\phi^3_{t=0}\log(\phi_{t=0})}\label{eq01:B_sol}.
\end{align}
where we have taken terms up to leading order in $\Btwo$. Numerically solving the equations of motion \eqref{eq01:non_dimensionalised_eom} in \cref{fig01:conserved}, we see that the expansions for $\Bzero$ and $\Btwo$ remain roughly constant\footnote{The can is rolling without slipping and hence energy is conserved.}. In the frictionless case, given in \Cref{app:slipping_rocking_can}, the two conserved quantities emerge naturally from the Lagrangian formulation of the equations of motion and correspond to components of the can's angular momentum. In the presence of friction this is no longer exactly the case, but similarities are explored in \cref{app:bzero_and_btwo}. 

\subsection{Reduced equation of motion}

If we substitute \cref{eq01:series_sols} into \eqref{eq01:Phi}, we obtain a second order nonlinear ODE in $\phi$ 
\begin{align}
	\ddphi =& \frac{a_3}{\phi^3} + \frac{a_{l2}\log\phi}{\phi^2} + \frac{a_{2}}{\phi^2} + \frac{a_{l1}\log\phi}{\phi} + \frac{a_{1}}{\phi} + a_{ll}(\log\phi)^2 + a_{l0}\log\phi + a_0  - 1 + \order{\phi(\log\phi)} \label{eq01:Phi_reduced_intermediate}
\end{align}
where trigonometric terms have been expanded in $\phi$. The coefficients $a_{ij}$ contain only terms in $\Bzero^2$, $\Btwo^2$ and $\Bzero\Btwo$ and depend on material parameters $h$, $k$, $a_p$ and $c_p$. We have also rescaled time by $\sqrt{a_p+1}$ and overloaded the notation.  The leading order coefficient $a_3= \Btwo^2a_p$ is positive. All coefficients are given in \cref{app:coefficients}.

In \cref{fig01:conserved}, $\Bzero$ and $\Btwo$ are small and of approximately the same order. Therefore, we assume
\begin{align}
  |\Btwo| = \epsilon^{1/2},\\
  \Bzero = \zeta\epsilon^{1/2} \label{eq01:epsilon_zeta}
\end{align}
where $\zeta = \order{1}$ and $0<\epsilon \ll 1$. Then we truncate \cref{eq01:Phi_reduced_intermediate} to find
\begin{align}
	\ddphi =& \epsilon\left(\frac{a_3}{\phi^3} + \frac{a_{l2}\log\phi}{\phi^2} + \frac{a_{2}}{\phi^2} + \frac{a_{l1}\log\phi}{\phi} + \frac{a_{1}}{\phi} + a_{ll}(\log\phi)^2 + a_{l0}\log\phi + a_0\right)  - 1,  \label{eq01:Phi_reduced}
\end{align}
where we have extracted a factor of $\epsilon$ from each relabelled coefficient. \Cref{eq01:Phi_reduced} can be written in Hamiltonian form; see  \cref{app01:hamiltonian}.

\begin{figure}
    \centering
    \begin{overpic}[percent,width=\textwidth]{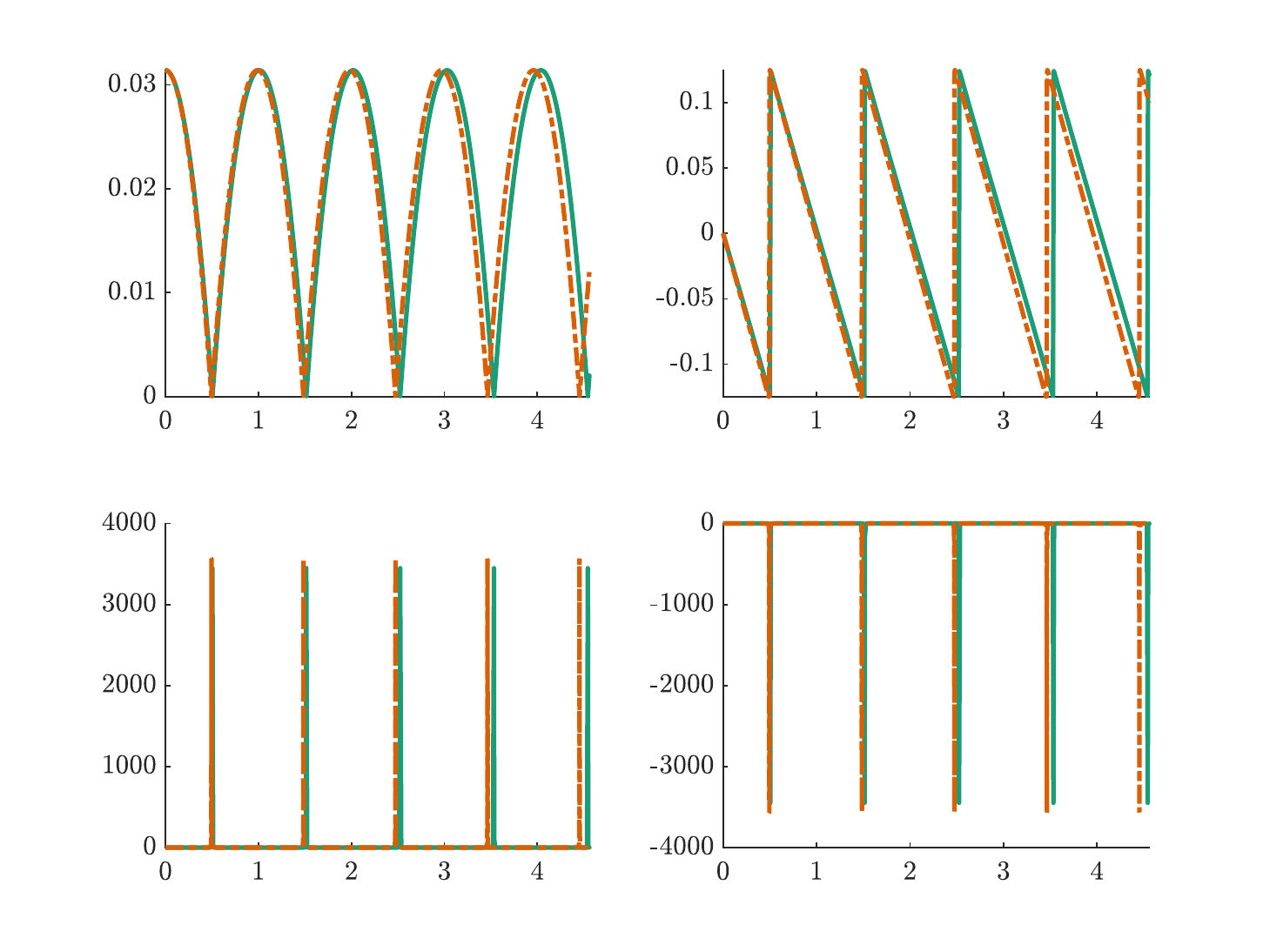}
    		\put(27,3){$t$}
    		\put(73,3){$t$}
    		\put(27,37){$t$}
    		\put(73,37){$t$}
    		\put(5,22){$\Psi$}
    		\put(5,55){$\phi$}
    		\put(50,22){$\Theta$}
    		\put(50,55){$\Phi$}
    \end{overpic}
    \caption{Numerical comparison of exact equations \cref{eq01:non_dimensionalised_eom} (solid green) and reduced equation \cref{eq01:Phi_reduced}. 
    }
    \label{fig01:eq_comparison}
\end{figure}

In Figure \ref{fig01:eq_comparison}, we compare the numerical solutions of the exact equations \cref{eq01:non_dimensionalised_eom} and the reduced equation \cref{eq01:Phi_reduced}. The initial conditions \cref{eq01:initial_conditions} and the scalings fix the values of $\epsilon =1.43\times 10^{-10}$ and $\zeta = 0.40$. The reduced equation \cref{eq01:Phi_reduced} shows some drift in the period, because the truncation removes the equilibrium $\phi= \phi^*$, so the can is unable to overturn. These errors are unimportant for the angle of turn phenomenon, which occurs when $\phi \ll \phi^*$. 

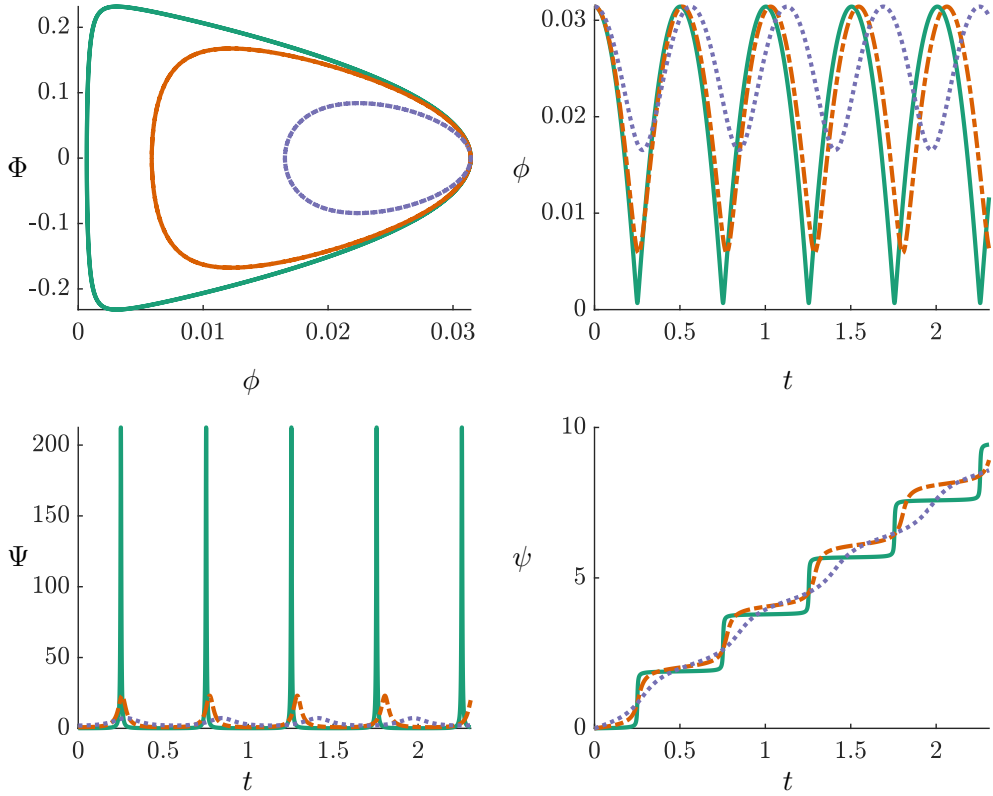
\begin{figure}
	\centering
	\begin{overpic}[percent,width=\textwidth]{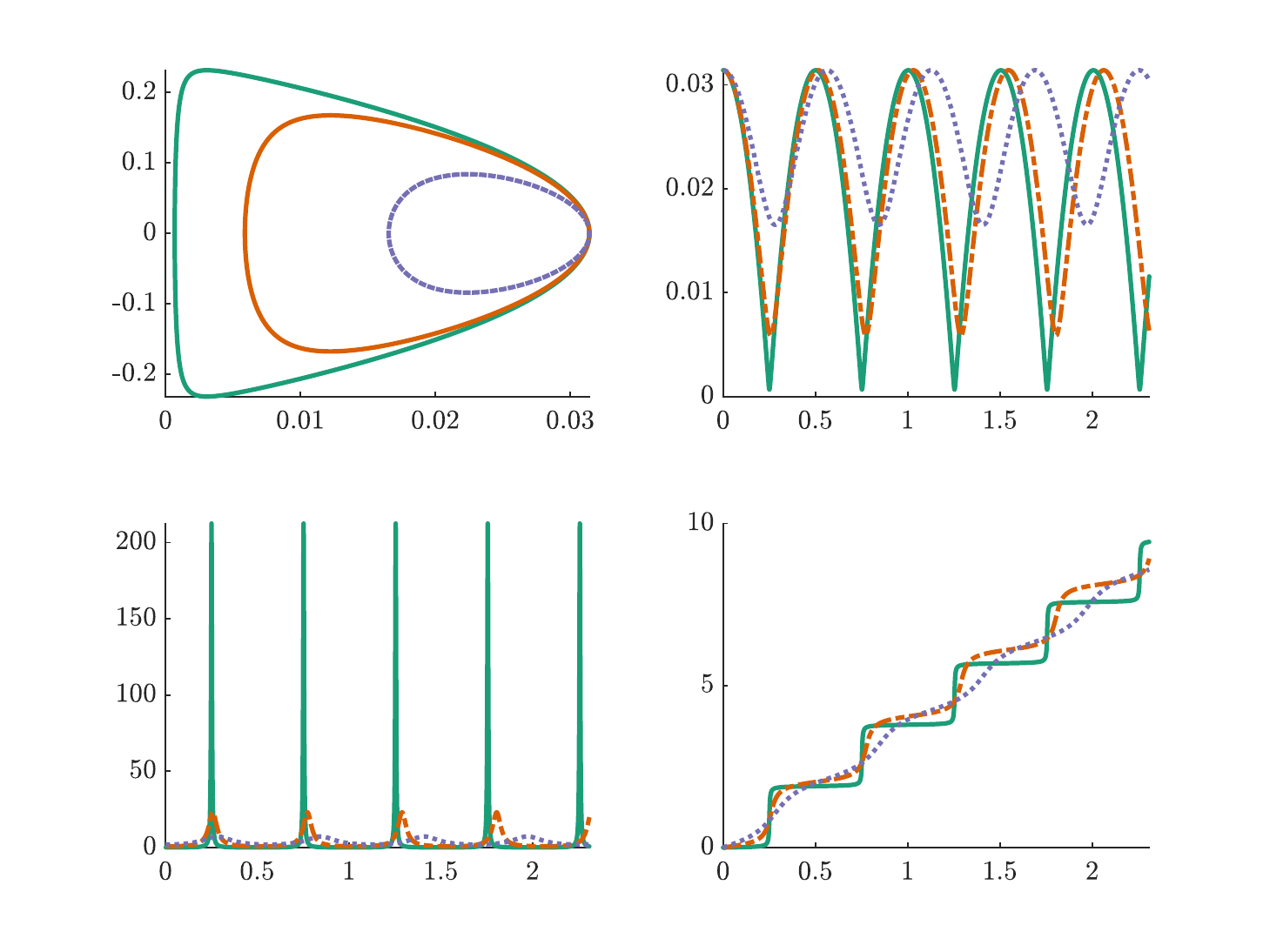}
			\put(27,3){$t$}
    		\put(73,3){$t$}
    		\put(27,37){$\phi$}
    		\put(73,37){$t$}
    		\put(7,22){$\Psi$}
    		\put(7,55){$\Phi$}
    		\put(50,22){$\psi$}
    		\put(50,55){$\phi$}
	\end{overpic}
    \caption{Time series and trajectories for \cref{eq01:Phi_reduced} with $\epsilon = 10^{-8}$ (green solid), $6.4\times10^{-7}$(orange dashed) and $4\times10^{-6}$(blue dotted) and $\zeta =1$. 
    }\label{fig01:sing_perturb}
\end{figure}

\Cref{fig01:sing_perturb} shows solutions of the reduced equation \cref{eq01:Phi_reduced} for different values of $\epsilon$, where we have used \cref{eq01:psi_series} to find $\Psi(\phi)$. 
As $\epsilon$ decreases, the nutation angle $\phi$ exhibits the expected bouncing behaviour,  $\psi$ shows step-like increases that correspond to the angle of turn, and the phase portrait shows clearly the repulsion from the singular line $\phi=0$ due to the $\phi^{-3}$ term. Thus $\epsilon$ mediates the angle of turn phenomenon.

\section{Asymptotic analysis}\label{sec01:asymptotic_analysis}

In this section we study the reduced equation \cref{eq01:Phi_reduced} using matched asymptotic expansions \cite{hinchperturbation}. 
When $\epsilon$ is small, \cref{fig01:sing_perturb} shows two regions of behaviour: when $\phi \gg 1$, the can acts like a compound pendulum (the outer region), and when $\phi \ll 1$, the bounce phenomenon occurs (the inner region). Matching the two solutions will yield a uniformly valid solution for one half period of the motion. 

\subsection{The outer solution}

When $\phi \gg 1$, we assume a regular perturbation in $\epsilon$
\begin{align}
    \phi(t) \sim \phi_0(t) + \epsilon\phi_1(t) + \order{\epsilon^2}.\label{eq01:outer_series_ansatz}
\end{align}
Upon substitution into the reduced equation \cref{eq01:Phi_reduced}, the leading order problem is $\ddphi_0 = -1$, which describes the can falling under gravity. 
If we take initial conditions $\phi(0) = I/2>0, \dphi(0) = 0$, we find
\begin{align}
    \phi_0(t) = \frac{I-t^2}{2}, \label{eq01:leading_order_outer}
\end{align}
which gives the expected quadratic form. So, in the absence of angular momentum in the $z^\mathcal{B}$ and $z^\mathcal{G}$ directions ($\epsilon=0$), the can falls flat at $t=\sqrt{I}$.

At $\order{\epsilon}$, using \cref{eq01:leading_order_outer}, we obtain
\begin{align}
	\ddphi_1 =& \frac{b_{3}}{(1 -\tau^2)^3} + \frac{b_{l2}\log(1 -\tau^2)}{(1 -\tau^2)^2} + \frac{b_{2}}{(1 -\tau^2)^2} + \frac{b_{l1}\log(1 -\tau^2)}{1 -\tau^2} + \frac{b_1}{(1 -\tau^2)} \label{eq01:order1_outer}\\
	& +b_{ll}\log^2(1 -\tau^2)  +  b_{l0}\log(1 -\tau^2) + b_0, \nonumber
\end{align}
with initial conditions $\phi_1(0) = \dphi_1(0) = 0$, and we have rescaled time $t = \sqrt{I}\tau$. The coefficients  in \cref{eq01:order1_outer} are given in \cref{app:b_coefficients}. 

\Cref{eq01:order1_outer} has the explicit solution
\begin{align}
	\phi_1(\tau) = \sum^3_{i=0}b_iJ_i(\tau)+  \sum^2_{n=0} b_{li}J_{li}(\tau) + b_{ll}J_{ll}(\tau), \label{eq01:outer_sol}
\end{align}
where 
\begin{subequations}\label{eq01:integrals}
\begin{align}
	J_i(\tau) =& \int\int  \frac{1}{(1-\tau^2)^i} \mathrm{d}\tau\mathrm{d}\tau, \quad (i=0 \ldots 3),\\
	J_{li}(\tau) = &  \int\int \frac{\log(1-\tau^2)}{(1-\tau^2)^i}\mathrm{d}\tau\mathrm{d}\tau, \quad (i=0 \ldots 2),\\
	J_{ll}(\tau) =& \int\int \log^2(1-\tau^2)\mathrm{d}\tau\mathrm{d}\tau.
\end{align}
\end{subequations}
Each double integral is evaluated in \cref{secA:integrals}. 
Therefore the outer solution up to and including $\order{\epsilon}$ terms is given by
\begin{align}
	\phi(t) \sim  \frac{I -t^2}{2} + \epsilon\left(\sum^3_{n=0}b_iJ_i\left(\frac{t}{\sqrt{I}}\right)+  \sum^2_{n=0} b_{li}J_{li}\left(\frac{t}{\sqrt{I}}\right) + b_{ll}J_{ll}\left(\frac{t}{\sqrt{I}}\right)\right) + \order{\epsilon^2}. \label{eq01:outer_solution}
\end{align}
As $t\to\sqrt{I}$, terms in $J_1(\frac{t}{\sqrt{I}}), J_3(\frac{t}{\sqrt{I}})$ become singular, so the outer solution is not valid as the can falls flat.

\subsection{The inner solution}

 Since the outer solution \cref{eq01:outer_series_ansatz} is not valid over the whole time domain, we require an inner solution. Let us define inner variables $\varphi, T$ as
 \begin{align}
     \phi = \epsilon^{1/2}\varphi, \quad t = \sqrt{I} + \epsilon^{1/2}T. \label{eq01:inner_scalings}
 \end{align}
 Substituting \cref{eq01:inner_scalings} into the reduced equation \eqref{eq01:Phi_reduced} and simplifying gives the inner problem 
 \begin{align}
	\varphi^{\prime\prime} =& \frac{a_p}{\varphi^3} + \epsilon^{1/2}\frac{a_{l2}\log\epsilon^{1/2}\varphi}{\varphi^2} + \epsilon^{1/2}\frac{a_{2}}{\varphi^2} + \epsilon\frac{a_{l1}\log\epsilon^{1/2}\varphi}{\varphi} + \epsilon\frac{a_{1}}{\varphi} + \epsilon^{3/2}a_{ll}\left(\log\epsilon^{1/2}\varphi\right)^2  \label{eq01:inner_full}\\
	&  + \epsilon^{3/2}a_{l0}\log\epsilon^{1/2}\varphi + \epsilon^{3/2}a_0 - \epsilon^{1/2}.\nonumber
\end{align}
where $\prime$ now denotes differentiation with respect the $T$. Taking $\varphi \sim \varphi_0 + \epsilon^{1/2}\varphi_1 + \order{\epsilon}$, we find to leading order 
 \begin{align}
     \varphi_0^{\prime\prime} = \frac{a_{p}}{\varphi^3_0}. \label{eq01:inner_problem}
 \end{align}
 \Cref{eq01:inner_problem} is equivalent to \cite[eq. (13)]{Srinivasan2008} and has solution
 \begin{align}
     \varphi(T) = \sqrt{P^2 + \frac{\left(\sqrt{a_{p}}T - Q\right)^2}{P^2}} + \order{\epsilon^{1/2}\log\epsilon},\label{eq01:inner_solution}
 \end{align}
 where $P$ and $Q$ are integration constants. Since both initial conditions have been used to find the integration constants for the outer solution, $P$ and $Q$ must be determined by matching.

\subsection{Matching}

In this section we match the outer and inner solutions to obtain a uniformly valid solution, using Van Dyke's matching rule \cite{hinchperturbation}. The inner solution \cref{eq01:inner_solution} written in the outer variable $t$ becomes

 \begin{align}
	   \frac{\sqrt{a_{p}}(\sqrt{I}-t)}{P} + \epsilon^{1/2}\frac{Q}{P} + \order{\epsilon}.
\end{align}
The outer solution \cref{eq01:outer_solution} written in the inner variable $T$ \cref{eq01:inner_scalings} becomes
 \begin{align}
	  -\sqrt{I}{}(t - \sqrt{I}) - \epsilon\frac{a_{p}}{2I^{5/2}(t - \sqrt{I})} +\order{\epsilon^{3/2}}.
\end{align}
Matching gives $P = \sqrt{a_p/I},\quad Q = 0$. Therefore the one term, matched inner solution is 
 \begin{align}
 	\varphi(T) = \sqrt{\frac{a_{p}}{I} + IT^2}.\label{eq01:phi_inner}
 \end{align}

\begin{figure}
    \centering
	\begin{overpic}[percent,width=\textwidth]{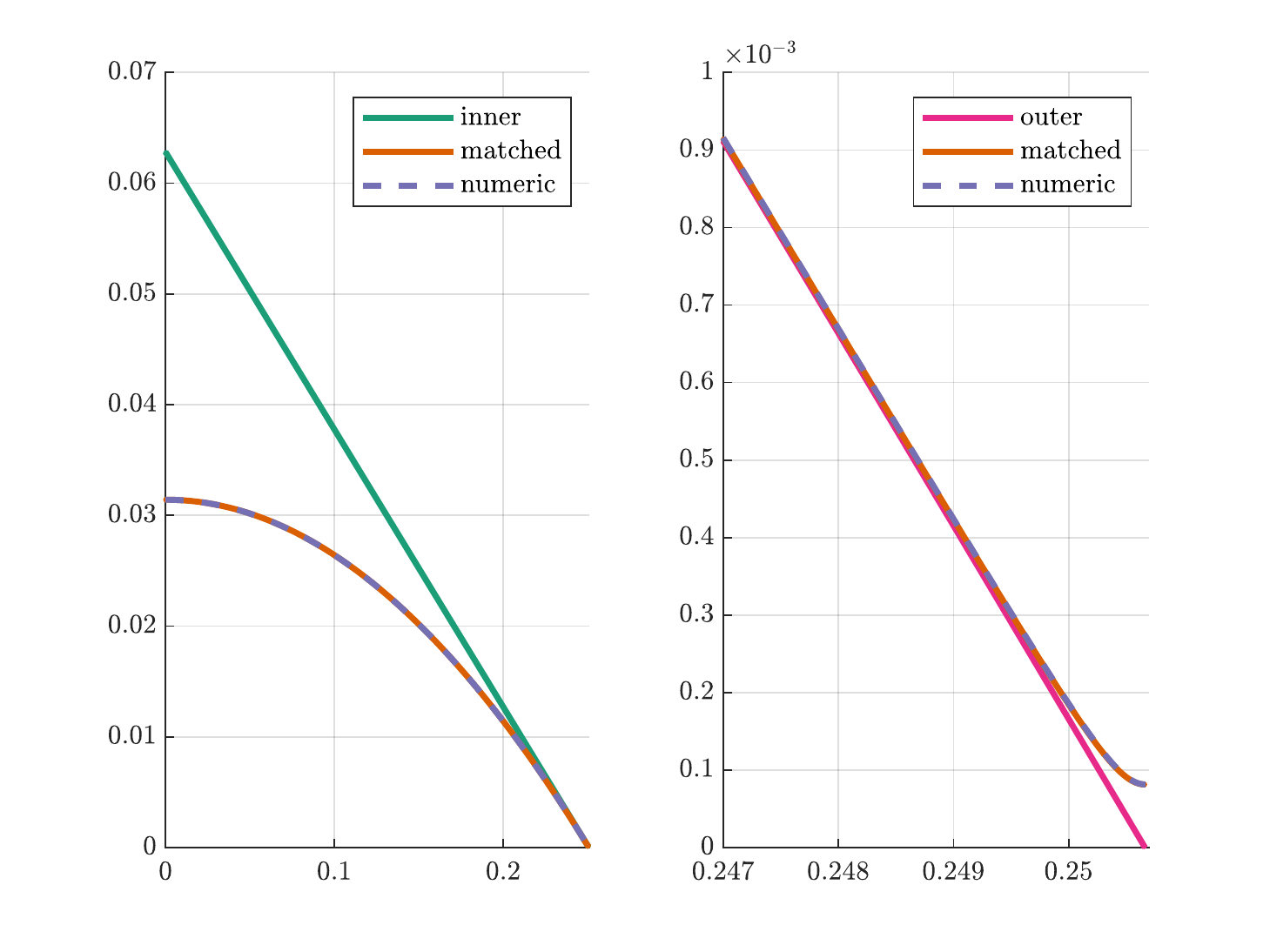}
			\put(30,3){$t$}
			\put(75,3){$t$}
			\put(4,40){$\phi$}
			\put(50,40){$\phi$}
	\end{overpic}
    \caption{Comparison of inner \eqref{eq01:inner_solution}, outer \eqref{eq01:outer_solution}, matched \eqref{eq01:uniformly_valid_solution} and numerical solution of the reduced equation \eqref{eq01:Phi_reduced}. Left: the time series for $0<t < \sqrt{I} = 0.2507$. The matched and numeric solutions lie on top of one another. Right: the time series for $0.247\textbf{}<t<0.2507$ showing the divergence of the outer solution from the matched and numeric solutions. 
    We take $\epsilon = \num{1.43e-10}$ and $I = \pi/50$.}
\label{fig01:inner_outer_comp}
\end{figure}
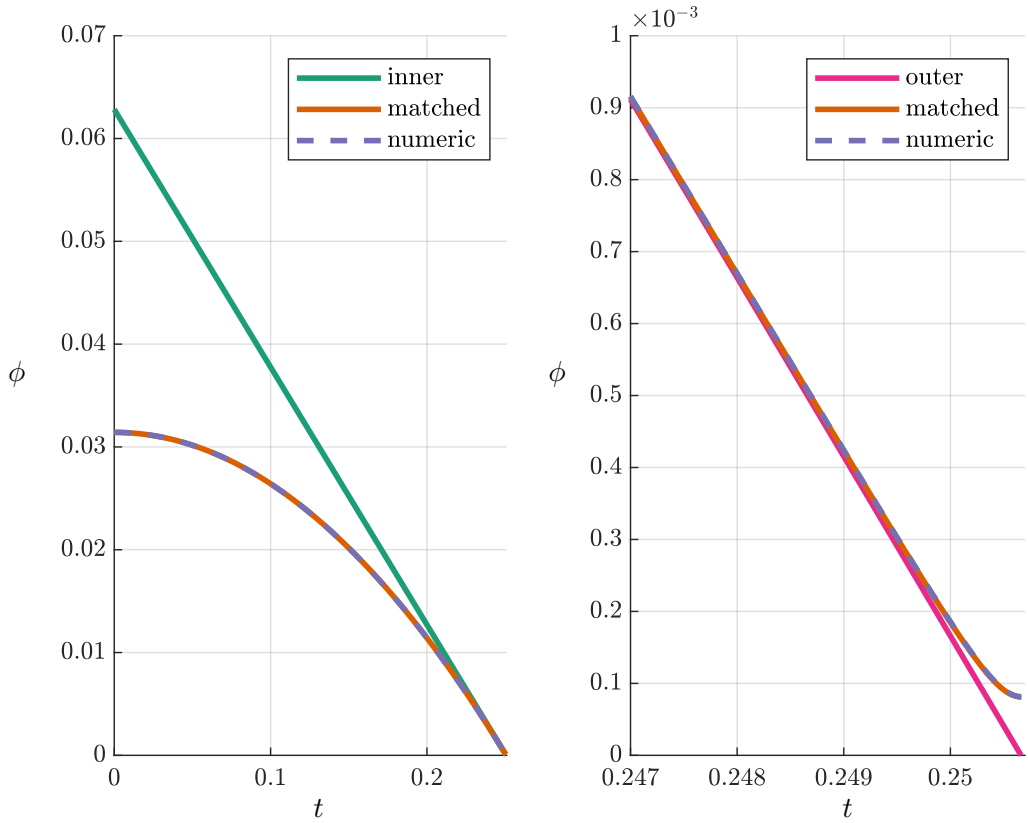

Hence the uniformly valid solution to the reduced equation \cref{eq01:Phi_reduced} can be written in the form
\begin{align}
	\phi(t)=&  -\frac{\left(\sqrt{I}-t\right)^2}{2}  + \sqrt{\frac{\epsilon a_{p}}{I} + I(t-\sqrt{I})^2} + \order{\epsilon}.\label{eq01:uniformly_valid_solution}
\end{align}
The inner, outer and matched solutions are shown in \cref{fig01:inner_outer_comp}, together with the numerical solution of the reduced equation \cref{eq01:Phi_reduced} for $\epsilon = \num{1.43e-10}$. The initial condition $\phi(0)$ determines $I$.  The matched solution is in excellent agreement with the numerical solution over the entire range of $\phi$, with the inner and outer solutions agreeing tangentially in their regions of applicability.

\subsection{Reconstructing the state variables}\label{sec:reconstruct_inner}

In the previous section we obtained a matched solution for the nutation angle $\phi(t)$ \cref{eq01:uniformly_valid_solution}. In this section, we reconstruct the remaining state variables of the can in the inner region, where the bounce takes place. 

We compute the inner solution for $\Psi$, using \cref{eq01:psi_series}. Recall that we had scaled time by a factor of $\sqrt{a_p+1}$ to obtain our reduced equation \cref{eq01:Phi_reduced}. We then shift to the inner variables $\phi = \epsilon^{1/2}\varphi$, $T = \epsilon^{-1/2}(t - \sqrt{I})$ and note that  $|\Btwo|=\epsilon^{1/2}$ to find
 \begin{align}
     \frac{1}{\sqrt{a_p+1}}\diff[T]{\psi} = \frac{\sign(\Btwo)}{\varphi^2} - \frac{\epsilon^{1/2}\sign(\Btwo)kh}{\varphi} + \order{\log\epsilon},
 \end{align}
Then, to leading order using \cref{eq01:phi_inner},
 \begin{align}
     \diff[T]{\psi} = \frac{\sqrt{a_p+1}\sign(\Btwo)}{\left(\frac{a_{3}}{I} + IT^2 \right)}. \label{eq01:Psi_inner}
 \end{align}
The maximum value of \cref{eq01:Psi_inner} occurs at $T=0$. Hence, reverting to the scalings in \cref{eq01:Phi_reduced} 
\begin{align}
     \diff[t]{\psi}\bigg|_\text{max} = \frac{I\sqrt{a_p+1}\sign(\Btwo)}{\epsilon^{1/2}a_p}. \label{eq01:dpsi_max}
 \end{align}
The presence of $\epsilon^{1/2}$ causes the large spikes visible in \cref{fig01:eq_comparison}, where $\diff[t]{\psi}\big|_\text{max}=3575$. \Cref{eq01:dpsi_max} gives $\diff[t]{\psi}\big|_\text{max} = 3449$, an error of 3.5\%. 
The same procedure is applied to the $\Theta$ expression \cref{eq01:theta_series} to give
  \begin{align}
     \diff[T]{\theta} = -\frac{\sqrt{a_p+1}\sign(\Btwo)}{\left(\frac{a_{3}}{I} + IT^2 \right)} = - \diff[T]{\psi}. \label{eq01:Theta_inner}
 \end{align}
From \cref{eq01:Psi_inner}, we find
\begin{align}
    \psi(T) &= \int \frac{\sqrt{a_p+1}\sign(\Btwo)}{\frac{a_{3}}{I} + IT^2}\mathrm{d}T  = \frac{\sqrt{a_p+1}\sign(\Btwo)}{a_p}\arctan\left(\frac{IT}{\sqrt{a_p}}\right) + C_\psi, \label{eq01:psi_inner}
\end{align}
and from \cref{eq01:Theta_inner}
\begin{align}
          \theta(T) &= - \frac{\sqrt{a_p+1}\sign(\Btwo)}{a_p}\arctan\left(\frac{IT}{\sqrt{a_p}}\right) + C_\theta, \label{eq01:theta_inner}
\end{align} 
where $C_\psi$ and $C_\theta$ are integration constants. From \cref{eq01:inner_solution}
 \begin{align}
     \diff[T]{\varphi} = \frac{IT}{\sqrt{\frac{a_{3}}{I} + IT^2}},
 \end{align}
thus completing the description of the  motion of the can during the bounce. Angles $\psi$ and $\theta$ exhibit a fast change over the bounce, corresponding to smoothed step functions, with the size of the step equal to the angle of turn $\Delta\psi$. The scaled angle $\varphi$ evolves according to a smoothed modulus function, avoiding $\varphi=0$. In unscaled variables the minimum angle reached by $\phi$ is $\phi_\text{min} = \sqrt{\epsilon a_{p}/I}$. 

\section{Physical Phenomena}\label{sec01:physical_phenomena}

In this section we use the expressions for the state variables in the inner region to understand the 
angle of turn phenomenon and its sign, the motion of the contact point and the coefficient of friction required for rolling without slipping to be sustained.

\subsection{Angle of turn}\label{sec01:angle_of_turn}

We compute $\Delta\psi$ by looking at the change of the inner solution for $\psi$ \cref{eq01:psi_inner}. This solution is only valid for $T<0$, but we can appeal to the symmetry of \cref{eq01:Psi_inner} to find
\begin{align}
    \Delta\psi = 2(\psi(0) - \psi(-\infty)) = \pi \sqrt{\frac{a_p+1}{a_{p}}}\sign(\Btwo). \label{eq01:angle_of_turn}
\end{align}
Reverting to unscaled units, the size of the angle of turn is
\begin{align}
    \left|\Delta\psi\right| = \pi\sqrt{\frac{A + mH^2 + mR^2}{A + mH^2}},\label{eq01:size_angle_of_turn}
\end{align}
agreeing with \cite{Srinivasan2008}, which used formal assumptions on $\psi$, $\theta$ and $\phi$. For our can, we calculate $|\Delta\psi| = 209^\circ$. The step-like behaviour of $\psi$ in the numerical solution, visible in \cref{fig01:non_dim_solution}, corresponds to multiple bounces, each with an angle of turn $|\Delta\psi| = 209^\circ$. The difference with respect to Srinivasan and Ruina is attributed to our choice of parameter values. 

From \cref{eq01:angle_of_turn}, $\lim_{a_p \to \infty}|\Delta\psi| \to \pi$, so rotation by any can is at least $\pi$ before rising back up, regardless of the material parameters. Conversely, $\lim_{a_p \to 0}|\Delta\psi|\to \infty$, suggesting that the can completes many revolutions before rising back up. The feasibility of such a large angle of turn is discussed in \cref{sec01:friction}.

\Cref{eq01:angle_of_turn} also contains information about the direction of the angle of turn determined by $\sign(\Btwo)$ \cite{Cushman2006}. We rewrite the expression for $\Btwo$ \cref{eq01:B_sol} as
\begin{align}
\Btwo =  \left (- \frac{kc_p}{2}\Theta_{t=0} + \Psi_{t=0}(1-\frac{kc_p}{2})\right) \phi_{t=0}^2 + \order{\phi_{t=0}^3} \label{eq01:b_sol_later}
\end{align}
To leading order, the sign of $\Btwo$ is determined by whether the initial conditions lie to the left or the right of the line
 \begin{align}
\Theta_{t=0} = \left(\frac{2}{c_pk} - 1\right)\Psi_{t=0}. \label{eq01:boundary}
\end{align}
 
\begin{figure}
	\centering
	\begin{overpic}[percent,width=0.8\textwidth]{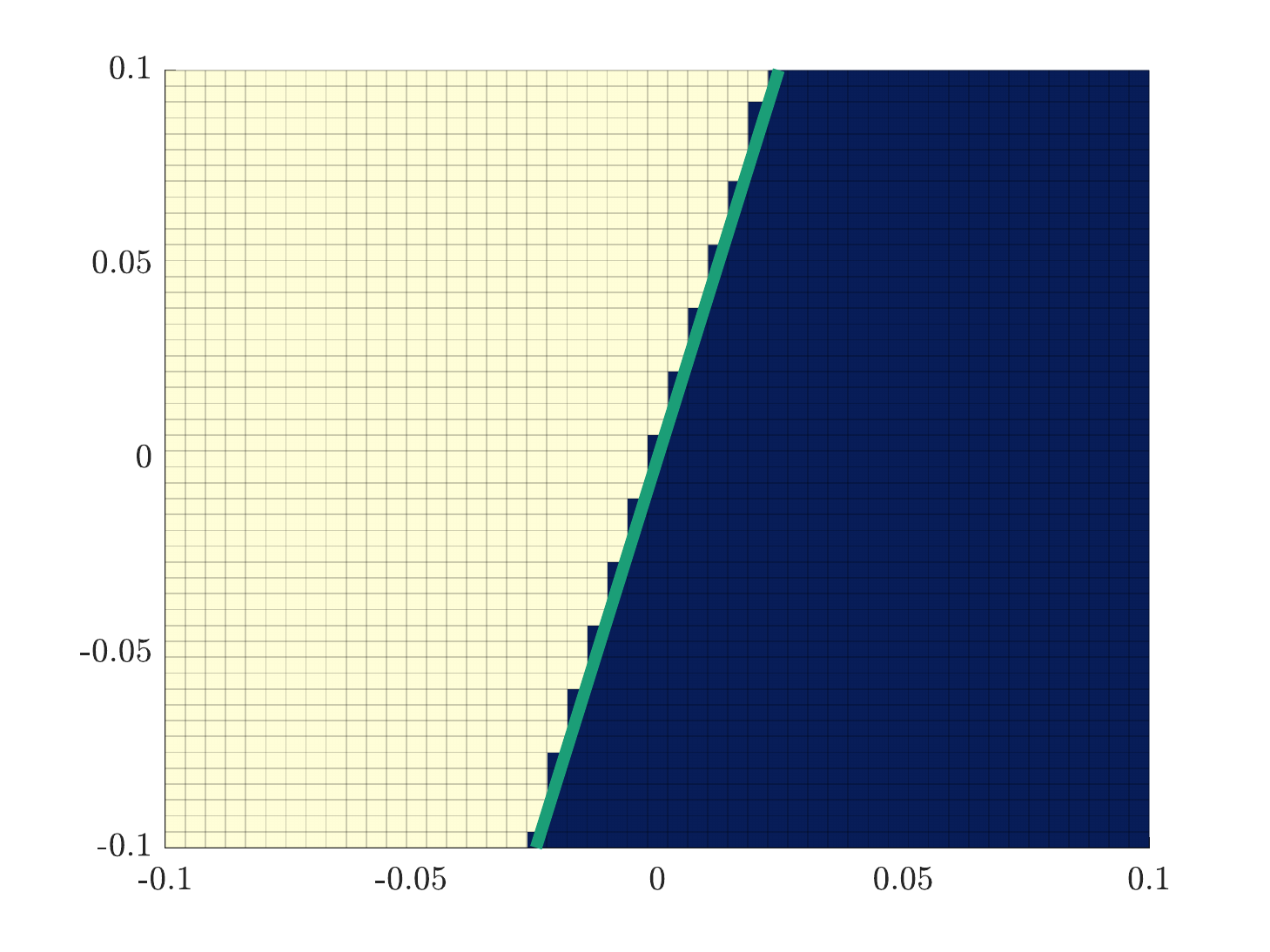}
		\put(0,40){$\dtheta$}
		\put(50,0){$\dpsi$}
	\end{overpic}
	\caption{Numerical simulations of the can with unscaled initial conditions $(\dpsi,\phi,\dphi,\dtheta)= (\dpsi_{t=0},\pi/100,0,\dtheta_{t=0})$. A square is coloured blue if $\Delta\psi <0$ and cream if $\Delta\psi > 0$. The lighter blue line is given by \cref{eq01:boundary}.}\label{fig01:left_or_right}
\end{figure}

In \cref{fig01:left_or_right}, we show excellent agreement between \eqref{eq01:boundary} and the numerical solutions of \cref{eq01:non_dimensionalised_eom}, for different initial conditions $\Psi_{t=0}$ and $\Theta_{t=0}$. 

\subsection{Contact point motion}\label{sec01:contact_pt}

As the can enters the bounce phase, the instantaneous contact point races quickly around the rim of the can. In \cite{Cushman2006}, paths traced out by the contact point are found by numerical integration of the equations of motion. In this section, we investigate these paths analytically, using the inner solutions in \cref{sec:reconstruct_inner}.

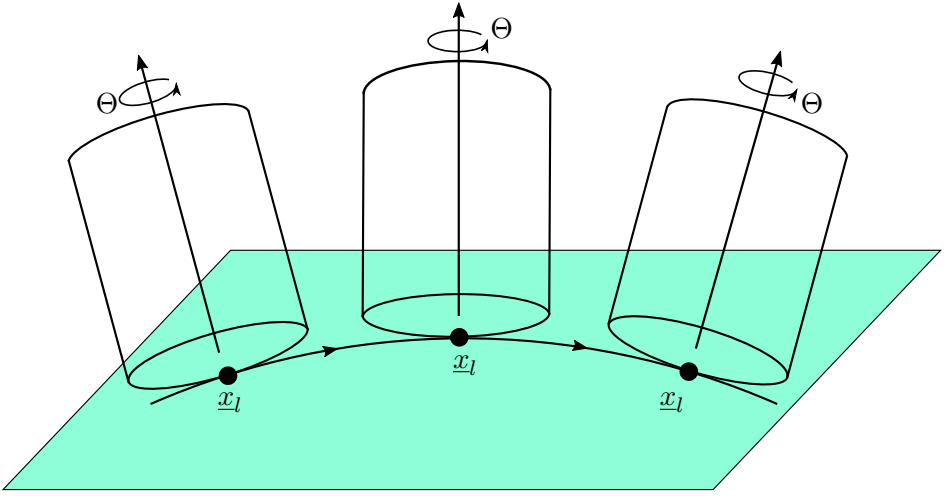
\begin{figure}
 	\centering
 	\begin{overpic}[width = 0.8\textwidth,percent]{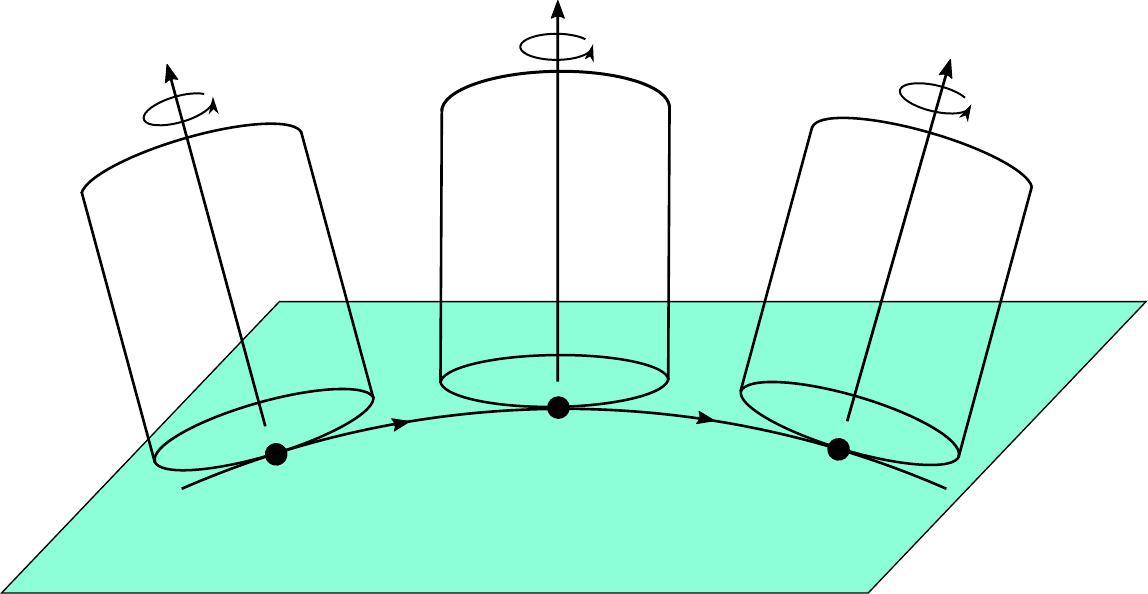}
 		\put(23,9){$\vec{x}_l$}
 		\put(48,13){$\vec{x}_l$}
 		\put(70,9){$\vec{x}_l$}
 		\put(10,40){$\Theta$}
 		\put(52,48){$\Theta$}
 		\put(85,40){$\Theta$}
 	\end{overpic}
 	\caption{
  The path that the instantaneous contact point traces out as the can moves over the plane is called the contact locus, $\vec{x}_l$.}\label{fig01:contact_locus}
 \end{figure}					 

We introduce the contact locus $\vec{x}_l$, the point of contact between the can and the plane (\cref{fig01:contact_locus}) located at $\vec{GP}^\mathcal{B}$ (\cref{fig01:can_set_up}).
$\vec{x}_l$ is not fixed to the can, but rotating with it. It moves in the plane with velocity $\vec{v}_{l}^\mathcal{G}$ given by 
\begin{align}
	\vec{v}_{l}^\mathcal{G} = \vec{v}_{G}^\mathcal{G}  +  \mat{R}_{\mathcal{B}\mathcal{G}}\left((\vec{\Omega}^\mathcal{B} - \Theta \hat{\vec{z}}^\mathcal{B})\times\vec{GP}^\mathcal{B}\right), \label{eq01:contact_locus}
\end{align}
The centre of mass velocity $\vec{v}_{G}^\mathcal{G}$ is given in \cref{eq01:com_velocity}. Applying scalings \cref{eq01:scalings} and disregarding the zero $z$ component, we find
\begin{align}
    \vec{v}_l^\mathcal{G} = \diff[t]{\vec{x}_l} = \Theta\begin{pmatrix}
        -\sin\psi.\\ \cos\psi
    \end{pmatrix}, \label{eq01:contact_point_velocity}
\end{align}
If we now convert to inner variables, noting that $\diff[T]{\psi} = - \diff[T]{\theta}$ from \cref{eq01:Theta_inner}, we have
\begin{align}
    \diff[T]{\vec{x}_l} = -\diff[T]{\psi}  \begin{pmatrix}
        -\sin\psi\\ \cos\psi
    \end{pmatrix},
\end{align}
which is integrable. Hence, to leading order, the inner solution for the contact locus describing the path traced out by the can is
\begin{subequations}\label{eq01:contact_pt_inner}
\begin{align}
  x_l(T) & = - \cos\left(\frac{\sqrt{a_p+1}\sign(\Btwo)}{a_p}\arctan\left(\frac{IT}{\sqrt{a_p}}\right) + C_\psi\right) + C_{xl},\\
  y_l(T) &= \sin\left(\frac{\sqrt{a_p+1}\sign(\Btwo)}{a_p}\arctan\left(\frac{IT}{\sqrt{a_p}}\right) + C_\psi\right) + C_{yl},
\end{align}
\end{subequations}
where $C_{xl}$ and $C_{yl}$ are constants of integration. Therefore, during the bounce phase the contact locus  moves in a circular arc, of size $|\Delta\psi|$, with the same radius as the can. 

To understand how the contact locus changes in the outer region we recall the series solutions for $\Theta$ and $\Psi$, \cref{eq01:psi_series,eq01:theta_series}. Rescaling time by $\sqrt{a_p+1}$, inserting the outer solution for $\phi$ \cref{eq01:outer_solution}, and discarding the $\order{\epsilon}$ terms we obtain
\begin{align}
	\diff{\psi} = \frac{4\epsilon^{1/2}\sign(\Btwo)\sqrt{a_p+1}}{(I-t^2)^2}  = -\diff{\theta}. \label{eq01:dpsi_outer_contact}
\end{align}
Integrating \cref{eq01:dpsi_outer_contact} gives the outer solution for $\psi$: 
\begin{align}
	\psi(t) = 2\epsilon^{1/2}\sign(\Btwo)\sqrt{a_p+1}\left( \frac{t}{I(I-t^2)} + \arctanh\left(\frac{t}{\sqrt{I}}\right)\right) + C_\psi,
\end{align}
and hence from \cref{eq01:contact_point_velocity}
\begin{subequations}\label{eq01:contact_pt_outer}
\begin{align}
  x_l(t) & = - \cos\left(2\epsilon^{1/2}\sign(\Btwo)\sqrt{a_p+1}\left( \frac{t}{I(I-t^2)} + \arctanh\left(\frac{t}{\sqrt{I}}\right)\right) + C_\psi\right) + C_{xl},\\
  y_l(t) &= \sin\left(2\epsilon^{1/2}\sign(\Btwo)\sqrt{a_p+1}\left( \frac{t}{I(I-t^2)} + \arctanh\left(\frac{t}{\sqrt{I}}\right)\right) + C_\psi\right) + C_{yl}.
\end{align}
\end{subequations}
The presence of $\epsilon^{1/2}$ in \cref{eq01:contact_pt_outer} shows that $\vec{x}_l$ barely changes in the outer solution until\footnote{In \cite{Cushman2006}, the contact locus of a thin disk ($H=0$) rotates clockwise at small $\phi$. After overturning, the disk rotates anticlockwise (or vice versa). This change of rotation direction is not visible in the inner and outer solutions \cref{eq01:contact_pt_inner,eq01:contact_pt_outer} because the thick disk cannot overturn.} $t \to \sqrt{I}$. 

To avoid multiple points of contact, $\phi\in (0,\pi/2)$. Both solutions \cref{eq01:contact_pt_inner,eq01:contact_pt_outer} rotate clockwise if $\sign(\Btwo)>0$ and anticlockwise if $\sign(\Btwo)<0$. 

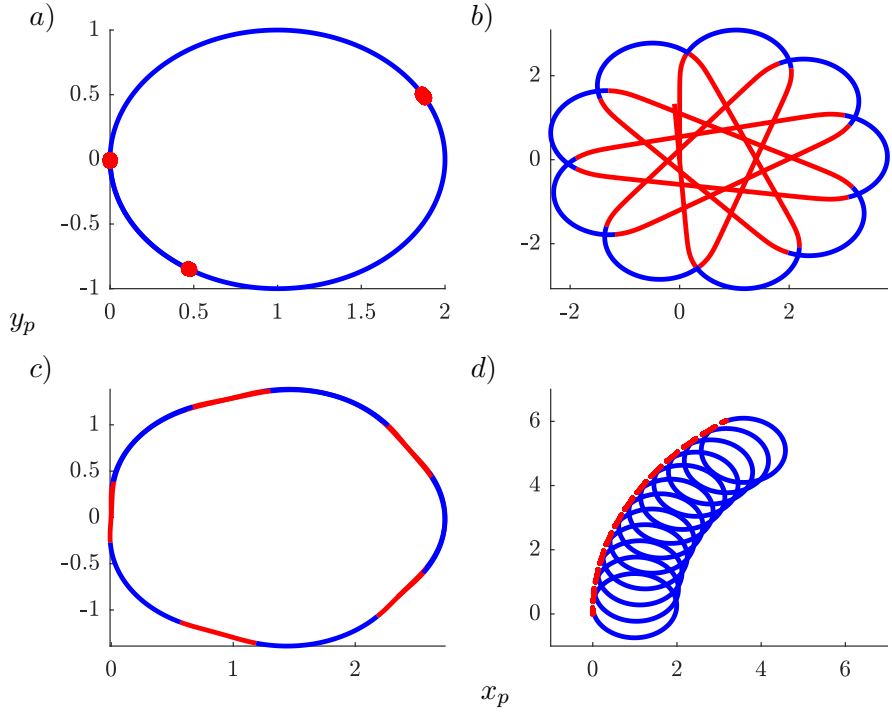
\begin{figure}
    \centering
    \begin{overpic}[scale=0.9,percent]{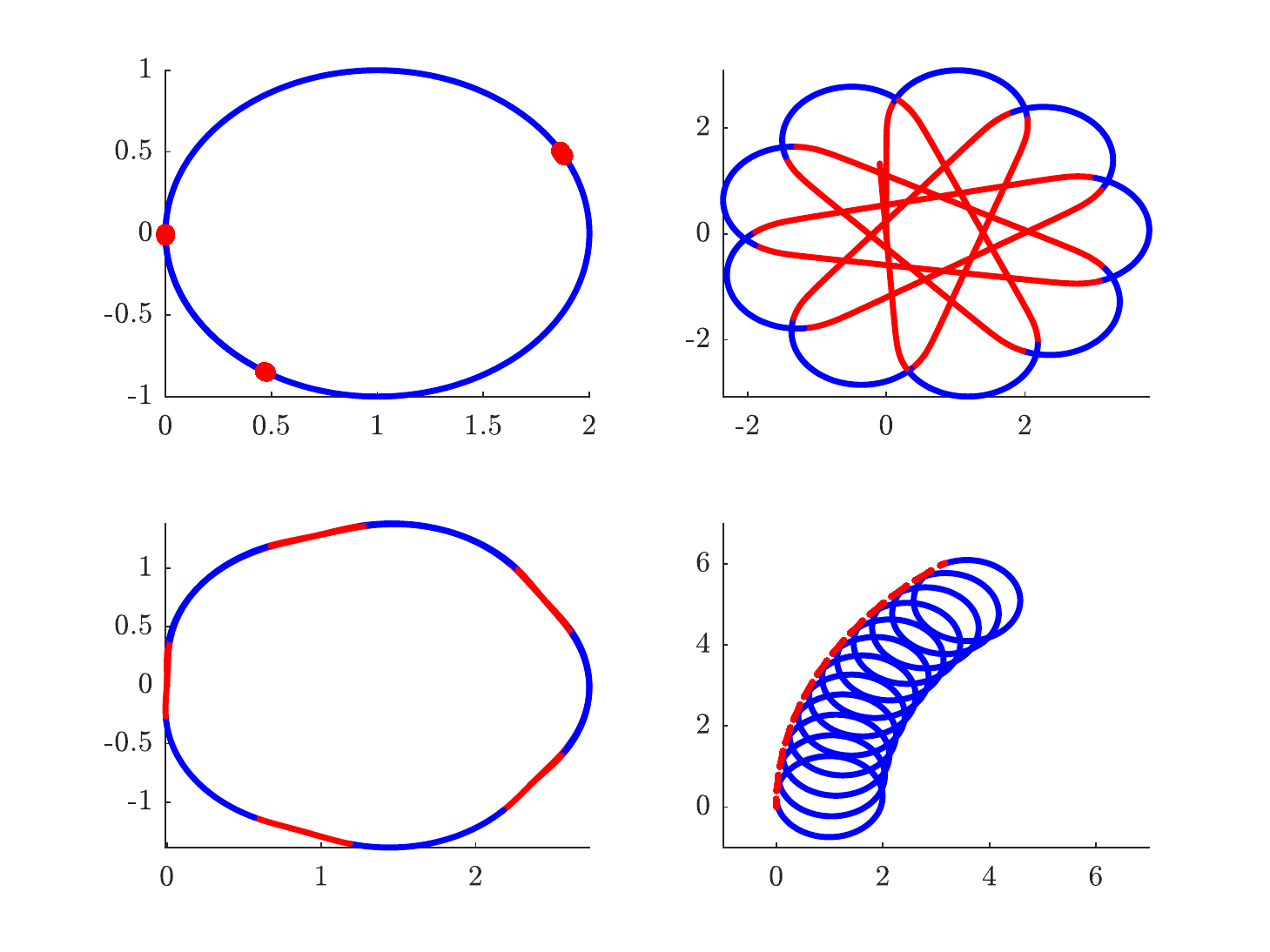}
        \put(50,3){$x_p$}
        \put(3,40){$y_p$}
        \put(5,70){$a)$}
        \put(49,70){$b)$}
        \put(5,35){$c)$}
        \put(49,35){$d)$}
    \end{overpic}
    \caption{Contact locus trajectories for different initial conditions. Portions of a trajectory where $\phi$ is small compared to the initial condition are coloured blue ($\phi<\phi_0/10$). To obtain the trajectories the contact point velocity \cref{eq01:contact_point_velocity} is integrated numerically using solutions of \cref{eq01:non_dimensionalised_eom}. The unscaled initial conditions are: a) as in\cref{eq01:initial_conditions}, b) $(\dpsi,\phi_0,\dphi,\dtheta) = (0,\phi^*-\num{1e-5},0,2)$, c)  $(\dpsi,\phi,\dphi,\dtheta) = (0.72667,\pi/100,0,7)$ and d) $(\dpsi,\phi,\dphi,\dtheta)= (1.7267,\pi/100,0,7)$. 
    }    \label{fig01:contact_pt_motion}
\end{figure}

\Cref{fig01:contact_pt_motion} shows the contact locus trajectory for different sets of initial conditions. The trajectories are calculated numerically using \cref{eq01:non_dimensionalised_eom}. In a), with small $\phi$ and small $\epsilon$, the contact locus trajectories are circular with large arcs in the inner solution and small arcs in the outer solution. Panel b) shows motion with initial conditions just off the balancing angle $\phi^* = \arctan(1/h)$ with a large $\dtheta$. The can rolls in a straight line before falling almost flat, turning around and repeating, producing a petaloid pattern. Such a pattern cannot be seen in the reduced system \eqref{eq01:Phi_reduced}, because it requires the presence of the saddle equilibrium at  $\phi^*$ which is destroyed by discarding $\order{\phi}$ terms. In c), we have small $\phi$,  large $\epsilon$. We see a roughly circular trajectory, but with wobbles caused by the large vertical angular momentum $\epsilon$. In d), the initial conditions are small $\phi$ with large $\dtheta$ and $\epsilon$. Close to the boundary defined by \cref{eq01:boundary}, we see circular movement, but along the top we see small cusp-like projections from the circle where the can reverses its direction. 

These trajectories have some striking patterns. Our analysis only applies to the  circular trajectory in a). 
Panel d) raises another question:  for what initial conditions  does the can's position remain bounded? If the initial conditions belong to the steady motion $S_\text{bal}$ \cref{eq01:equilibrium}, with $\Psi=0$, then the can rolls in a straight line with velocity $R\Theta$ and the position is unbounded.  In \cref{fig01:contact_locus} d), the contact locus appear to be proscribe a larger circle, but it is unclear if the position remains bounded. In the case of the thin disk, Borisov et al. \cite{Borisov2005} show that the contact locus is bounded for almost all initial conditions.

\subsection{Coefficient of friction}\label{sec01:friction}

To be valid, the equations of motion \cref{eq01:non_dimensionalised_eom} require that the coefficient of friction $\mu > |\vec{F(t)}|/N(t)$ for all time. Srinivasan and Ruina \cite{Srinivasan2008} found that the maximum value of $|\vec{F(t)}|/N(t)$ remains finite. In this section we investigate how  $|\vec{F}|/N$ depends on the material characteristics of the can.

\begin{figure}
	\centering
	\begin{overpic}[percent,width=0.8\textwidth]{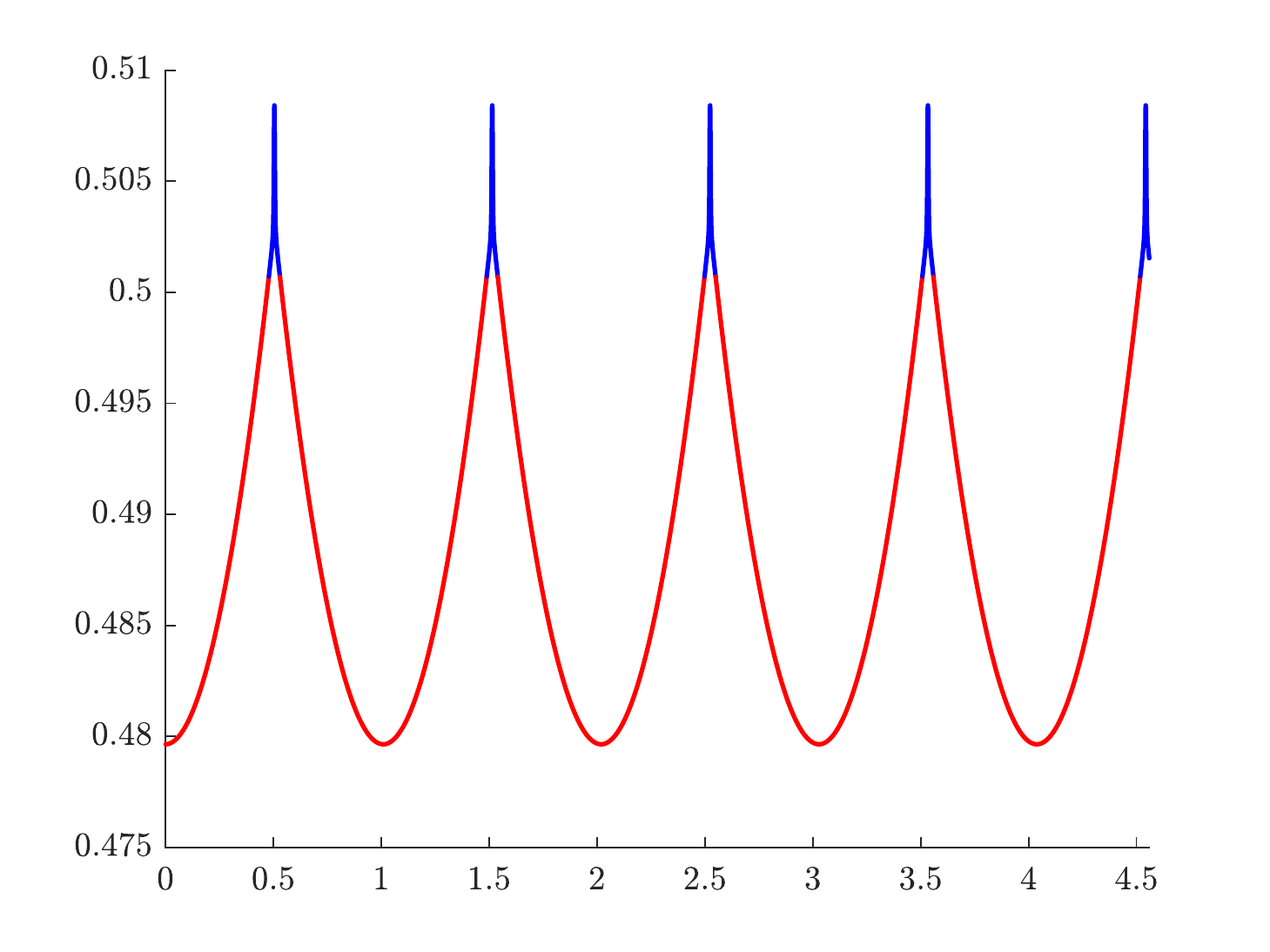}
			\put(50,0){$t$}
			\put(0,38){$|\vec{F}|/N$}
	\end{overpic}
	\caption{The ratio \cref{eq01:friction_short} versus time.  The plot is coloured blue if $\phi < \phi_0/10$ (the inner region), and red otherwise (the outer region). 
    }.
 \label{fig01:friction}
\end{figure}

The contact forces $F_x$, $F_y$ and $N$ are computed in \cref{eq01:fx,eq01:fy,eq01:normal} and 
\begin{align}
	\frac{|\vec{F}|}{N} = \frac{\sqrt{F_x^2 + F_y^2}}{N}. \label{eq01:friction_short}
\end{align}
 In \cref{fig01:friction}, we plot the friction ratio \eqref{eq01:friction_short} for a typical rocking can motion. The peak of the ratio $|\vec{F}|/N$ in \cref{fig01:friction} is $0.508$, suggesting that if $\mu > 0.508$ the can will roll without slipping.
 In this case, the outer region (in red, where $\phi=\order{1}$) requires less friction than the inner region (in blue, where $\phi \ll 1$)\footnote{This is not true for all initial conditions and parameter values. The outer region may require a higher coefficient of friction than the inner, particularly if the initial conditions have large $\Phi$.}. 

To look at the ratio of $\frac{|\vec{F}|}{N}$ over the bounce (inner region) we write \cref{eq01:friction_short} in inner variables and use \cref{eq01:phi_inner,eq01:Psi_inner,eq01:Theta_inner} along with  their derivatives. Expanding and keeping only the leading order terms in $\epsilon$ gives an approximate lower bound on the required coefficient of friction, given  by  
\begin{align}
    \frac{|\vec{F}|}{N} \approx \frac{h}{a_p} = \frac{mRH}{A + mH^2}
     \label{eq01:required_friction}.
\end{align} 
 Unlike in \cref{fig01:friction}, \cref{eq01:required_friction} has no dependence on time or initial conditions, which appear at higher order. 
 \Cref{eq01:required_friction} predicts a required coefficient of friction of $0.508$, agreeing with the numerical value. 

In \cref{eq01:angle_of_turn}, the angle of turn $|\Delta\psi|$ can be arbitrarily large if $a_p$ is small. However, from \cref{eq01:required_friction}, a small $a_p$ results in a large coefficient of friction, rendering large angles of turn infeasible. This explains why large angles of turn are not seen in practice: the can will slip instead. 

For a uniform density cylinder, the required coefficient of friction is simply
\begin{align}
    \frac{|\vec{F}|}{N} \approx  \frac{RH}{\frac{1}{12}(3R^2 + H^2) + H^2}  = \frac{12h}{3 + 13h^2}.
\end{align} 
This has a maximum value of $\mu = 0.96$, when $h = \sqrt{39}/13 \approx 0.48$, which is not often achieved in tabletop experiments. 

\section{Conclusion}

We have considered the problem of a can rolling on a rough horizontal plane, with nutation angle $\phi$. We reduced the problem to a second order ODE \cref{eq01:frob_ode}, which has a regular singularity at $\phi=0$. We then found a Frobenius solution \cref{eq01:psi_orig} involving two coefficients $\Bzero$ and $\Btwo$, which are related to angular momenta about the global vertical and symmetry axes, justifying the formal assumptions made by Srinivasan et al. \cite{Srinivasan2008}. Setting $\Bzero=\Btwo=0$ leads to a singular perturbation problem\footnote{The singular limit ($\epsilon=0$) yields a flat falling solution, which is studied by Cushman and Duistermaat \cite{Cushman2006} for the thin disk.}, prompting the introduction of $\epsilon$, a small parameter describing a combination of the angular momenta.  

The rocking can exhibits two distinct phenomena: for $\phi=\order{1}$, behaviour very similar to an inverted pendulum, and for $\phi \ll 1$, dynamics with the angle of turn.

This distinction allows us to use matched asymptotic expansions, with an outer region $\phi=\order{1}$ and an inner region $\phi \ll 1$, to derive a uniformly valid solution \cref{eq01:uniformly_valid_solution} that is in excellent agreement with numerical calculations of the reduced system \cref{eq01:Phi_reduced}. The solution of the inner problem was used to investigate of the angle of turn phenomenon. We computed the minimum angle $\phi_\text{min}$ achieved by the can  and the maximum angular velocity $\dpsi_\text{max}$ attained over the bounce. We recomputed the angle of turn $|\Delta\psi|$ derived by \cite{Srinivasan2008} and gained more information about the direction of the angle of turn. 
These key characteristics of the dynamics are backed up by numerical solutions of the full nonlinear equations \cref{eq01:non_dimensionalised_eom} in \textsc{matlab}. We also examine the motion of the contact locus $\vec{x}_l$ and see a range of different trajectories, from circular to petaloid motion and even cusp-like behaviour. 

Finally, we used the solution to the inner problem to obtain an approximate lower bound for the required coefficient of friction to avoid slip. To leading order, the lower bound is independent of the initial conditions and dependent only on the material characteristics of the can. For a typical can, a coefficient of friction $\mu \approx 0.51$ is required to avoid slipping.


An interesting extension to the rocking can problem is the addition of a forced horizontal plane. In the planar case, Hogan \cite{Hogan1990} explored the dynamics of a rigid rectangular block rocking and impacting with a sinusoidally-forced horizontal plane. The system contains a range of dynamical behaviour including period-doubling cascades. Rigid cylinders in three dimensions, such as classical columns or grain silos, may experience the same types of behaviour and have been explored in the structural engineering literature \cite{Stefanou2010,Vassiliou2017}. Such work may benefit from the approach taken in this paper.

\appendix
\section{Coefficients in \texorpdfstring{\cref{eq01:theta_orig,eq01:psi_orig}}{(4.12) and (4.13)}}
\label{app:coefficients_Psi_Theta}
\begin{align*}
    \Psi_{00} &= B_0, \\
    \Psi_{l0} &= \frac{\Btwo((kh)^2 - k)}{2},\\
    \Psi_{l1} &= \frac{\Btwo kh((kh)^2 - k)}{6},\\
    \Theta_{00} &= \frac{2-kc_p}{kc_p}\Bzero + \frac{3k(c_p + 1) -4 - 9h^2k^2}{6kc_p}\Btwo\\
    \Theta_{l0} &= \frac{\Btwo(h^2k-1)}{2c_p}(c_pk-2),\\
    \Theta_{l1} &= \frac{\Btwo(h^2k-1)}{2c_p}\frac{kh(c_pk-6)}{3}.\\
\end{align*}

\section{The frictionless case} \label{app:slipping_rocking_can}

In the frictionless case, a Lagrangian approach is more appropriate due to the absence of the non-holonomic constraint on the contact velocity. With no lateral force we can assume that the horizontal component of the velocity of the centre of mass is zero. The Lagrangian, $\mathcal{L} = T-V$ is 
\begin{align}
	\mathcal{L} &= \frac{\mathbf{I}\vec{\Omega
}\cdot\vec{\Omega}}{2} + \frac{m\dot{z}^2}{2} - mg(R\sin\phi + H\cos\phi)\nonumber\\
&= \frac{A}{2}(\dpsi^2\sin^2\phi + \dphi^2) + \frac{C}{2}(\dpsi\cos\phi + \dtheta)^2\\
& + \frac{m}{2}(\dphi(R\cos\phi - H\sin\phi))^2 - mg(R\sin\phi + H\cos\phi)
\end{align}
The equations of motion are found using the Lagrangian equations for the generalised coordinates $q\in \{\phi,\psi,\theta\}$
\begin{align}
	\frac{\mathrm{d}}{\mathrm{d}t}\frac{\partial \mathcal{L}}{\partial \dot{q}} - \frac{\partial \mathcal{L}}{\partial q} = 0.
\end{align}
Since both $\psi$ and $\theta$ do not appear in the Lagrangian, we have
\begin{align}
\frac{\mathrm{d}}{\mathrm{d}t}\left(C(\dpsi\cos\phi + \dtheta)\right) = 0.\\
	\frac{\mathrm{d}}{\mathrm{d}t}\left(A\dpsi\sin^2\phi + C\cos\phi(\dpsi\cos\phi + \dtheta)\right) = 0
\end{align}
Integration yields the two conserved quantities $H_z^{\mathcal{B}}$ \cref{eq01:ang_mom_B} and $H_z^{\mathcal{G}}$ \cref{eq01:ang_mom_G}. $H_z^{\mathcal{B}}$ corresponds to angular momentum about the symmetry axis of the can, and $H_z^{\mathcal{G}}$ corresponds to a combination of the angular momenta about the symmetry and vertical axes.

The equation of motion for $\phi$ is
\begin{align}
	(A + m(R\cos\phi - H\sin\phi)^2)\ddphi& =  m\dphi^2(R\sin\phi + H\cos\phi)(R\cos\phi - H\sin\phi) \\ 
 +  \dpsi(A -C)\sin\phi\cos\phi - &C\dtheta\dpsi\sin\phi - mg(R\cos\phi - H\sin\phi).\nonumber
\end{align}
Applying the scalings from \eqref{eq01:scalings},
and eliminating $\dpsi, \dtheta$ through \cref{eq01:ang_mom_B,eq01:ang_mom_G}, gives a planar ODE.
Introducing, as in \cref{eq01:epsilon_zeta},  $H_z^{\mathcal{B}}=\epsilon^{1/2} \ll 1$,    $H_z^{\mathcal{G}} = \zeta \epsilon^{1/2}$ where $\zeta = \order{1}$, and expanding in $\epsilon$ gives, on truncation
\begin{align}
	\ddphi = \epsilon\left(\frac{\alpha_{3}}{\phi^{3}} +  \frac{\alpha_{2}}{\phi^{2}} + \frac{\alpha_{1}}{\phi} + \alpha_0\right)  + \frac{h\dphi^2 -1}{a +1},\label{app01:Phi_reduced_slip}
\end{align}
with
\begin{align}
	\alpha_{3} = \frac{(\zeta-1)^2}{a(a + 1)}.
\end{align}
\Cref{app01:Phi_reduced_slip} is the slipping equivalent of the reduced equation \cref{eq01:Phi_reduced}.

\section{Physical justification for $\Bzero$ and $\Btwo$}
\label{app:bzero_and_btwo}

Here we give a physical understanding of $\Bzero$ and $\Btwo$ in \cref{eq01:A_sol,eq01:B_sol}. Consider the angular momenta about the $z^{\mathcal{B}}$ and $z^{\mathcal{G}}$ axes, given respectively by
\begin{align}
	H_z^\mathcal{B} &= c(\Theta + \Psi\cos\phi) \approx c(\Theta + \Phi)  - \frac{c\Psi\phi^2}{2} + \order{\phi^{4}}, \label{eq01:ang_mom_B}\\
	H_z^\mathcal{G} &= a\Psi\sin^2\phi + c\cos\phi(\Psi\cos\phi + \Theta) \approx c(\Theta + \Psi) + \left(-\frac{c}{2}(\Theta+ 2\Psi) + a\Psi\right)\phi^2 + \order{\phi^4}. \label{eq01:ang_mom_G}
\end{align}

For small $\phi$, $\Bzero -  \frac{kc_p}{2c} H_z^\mathcal{B}  = \order{\phi_{t=0}}$, suggesting that $\Bzero$ is approximately proportional to the angular momentum about the body symmetry axis. Correspondence with angular momentum is less clear for $\Btwo$, it bears similarities to the quantity $H_z^\mathcal{B} - H_z^\mathcal{G}$. We suggest that the conserved quantities $\Bzero$ and $\Btwo$ correspond to angular momenta.

\begin{figure}
	\centering
	\begin{overpic}[percent,width=0.5\textwidth]{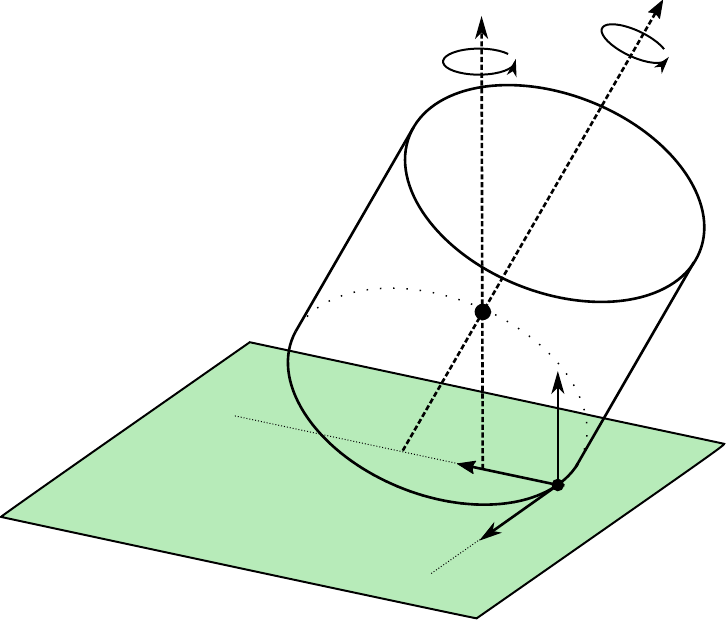}
			\put(66,6){$F_y^\mathcal{I}$}
			\put(58,18){$F_x^\mathcal{I}$}
			\put(73,32){$N$}
			\put(60,41){$\vec{G}$}
			\put(56,73){$\Psi$}
			\put(88,73){$\Theta$}
			\put(92,83){$z^\mathcal{B}$}
			\put(62,82){$z^\mathcal{G}$}
	\end{overpic}
	\caption{The contact forces on the can. The lines of action of both $N$ and $F_x^\mathcal{I}$ intersect with the axes $z^\mathcal{G}$ and $z^\mathcal{B}$ and, therefore, cannot affect the angular momenta about those axes. Only the tangential component of friction $F_y^\mathcal{I}$ can affect the angular momenta and this is small \cref{eq01:friction_tangential}.}\label{fig01:conserved_quantity_diagram}
\end{figure}
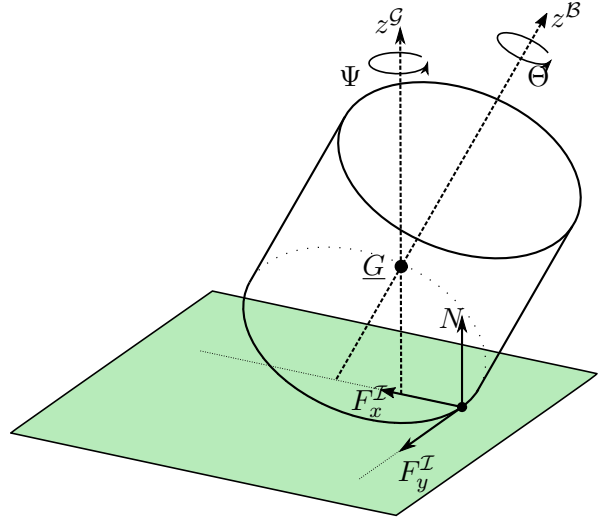

In the frictionless case, $H_z^\mathcal{B}$ and $H_z^\mathcal{G}$ emerge naturally as conserved quantities from the Lagrangian formulation of the equations of motion (see \cref{app:slipping_rocking_can}). Consider \cref{fig01:conserved_quantity_diagram}: without friction the only contact force is the normal force $\vec{N}$,  acting in the same plane as $H_z^\mathcal{B}$ and $H_z^\mathcal{G}$. Angular momentum must, therefore, be conserved in this case.

In the presence of friction, angular momentum is not conserved. \Cref{fig01:conserved_quantity_diagram} shows the component $F_y^\mathcal{I}\ne 0$, which acts in the tangential direction, affecting angular momentum, and given by
\begin{align}
    F_{y}^\mathcal{I} = (\mat{R}_{\mathcal{G}\mathcal{I}}\vec{F}^\mathcal{G})\cdot \hat{\vec{y}}^\mathcal{I} =  mgk\Phi(a\Psi\sin\phi + ch(\Psi\cos\phi + \Theta)). \label{eq01:friction_tangential}
\end{align}
For small $\phi$, $|F_{y}^\mathcal{I}| \propto (\Theta + \Psi$), which is in turn proportional to $\Bzero$, to leading order. If $\Bzero$ is also small, then the effects of the tangential force on the angular momenta are small.

\section{Coefficients in \texorpdfstring{\cref{eq01:Phi_reduced_intermediate}}{(4.21)}}
\label{app:coefficients}

\begin{align*}
	a_3 =& \Btwo^2a_p,\\
	a_2 =& -2\,{\frac {h \left(  \left( -1/3+ \left( -3/4\,{h}^{2}+{\it a_p}\,{
\it c_p} \right) {k}^{2}+ \left( -{\it c_p}/2+1/4 \right) k \right) \Btwo+\Bzero
 \right) \Btwo}{k{\it c_p}}},\\
 a_1 =& \frac{2}{kc_p}\Bigg[\bigg(\Big(\big( 1/3+1/2\,{h}^{2}{k}^{3}{\it a_p}-1/4
\,{k}^{2}{h}^{2}+ ( -{\it a_p}/3-1/4) k \big) \Btwo \\
& +\Bzero ( 
{\it a_p}\,k-1 )  \Big) {\it c_p}-\Btwo{h}^{2}{k}^{2} ( {h}^{2}
k-1)  \bigg) \Btwo\Bigg],\\
a_0 =& \frac {1}{2880\,k{\it c_p}} \Bigg[ \left( -125\,{\it c_p}\,{k}^{5}+
1440\,{k}^{4} \right) {\Btwo}^{2}{h}^{5}\\
&-120\,{k}^{2}B \left(  \left( 34-{
\frac {32\,{\it c_p}\,{k}^{2}}{3} \left( {\it a_p}+{\frac{11}{128}}
 \right) }+ \left( {\frac {293\,{\it c_p}}{12}}+36 \right) k \right) \Btwo+
\Bzero \left( k{\it c_p}-84 \right)  \right) {h}^{3}\\
&+\left(  \left[ -960+
 \Big(135 -3200\,{\it a_p} \right) {\it c_p}\,{k}^{3}+ \left( 5850+1280\,{
\it a_p} \right) {\it c_p}\,{k}^{2}+ \left(720 -1032\,{\it c_p}
 \right) k \right] B_{-2}^2\\
 &-3840\, \left[ -5/4+{\it c_p}\, \left( {\it 
a_p}+{\frac{3}{32}} \right) {k}^{2}+ \left( -{\frac {27\,{\it c_p}}{16}}
+3/8 \right) k \right] \Bzero\Btwo-5760\,{\Bzero}^{2} \Big) h\Bigg],\\ 
	a_{l2} =&{\frac {h{\Btwo}^{2} \left( {h}^{2}k-1 \right) }{{\it c_p}}},\\
	a_{l1} =& - \left( {h}^{2}k-1 \right) {\Btwo}^{2} \left( {\it a_p}\,k-1 \right),\\
	a_{l0} =&2\,{\frac {h \left( {h}^{2}k-1 \right) \Btwo}{{\it c_p}} \left(  \left( -{
\frac{5}{12}}+ \left( -3/4\,{h}^{2}+1/3\,{\it a_p}\,{\it c_p} \right) {k
}^{2}+ \left( -{\it c_p}/2+1/4 \right) k \right) \Btwo+\Bzero \right) },\\
a_{ll} =& -{\frac {hk \left( {h}^{2}k-1 \right) ^{2}{\Btwo}^{2}}{2\,{\it c_p}}}.
\end{align*}

\section{Hamiltonian system} \label{app01:hamiltonian}
\begin{figure}
\centering
\begin{minipage}{0.45\textwidth}
	\centering
	\begin{overpic}[percent,width=\textwidth]{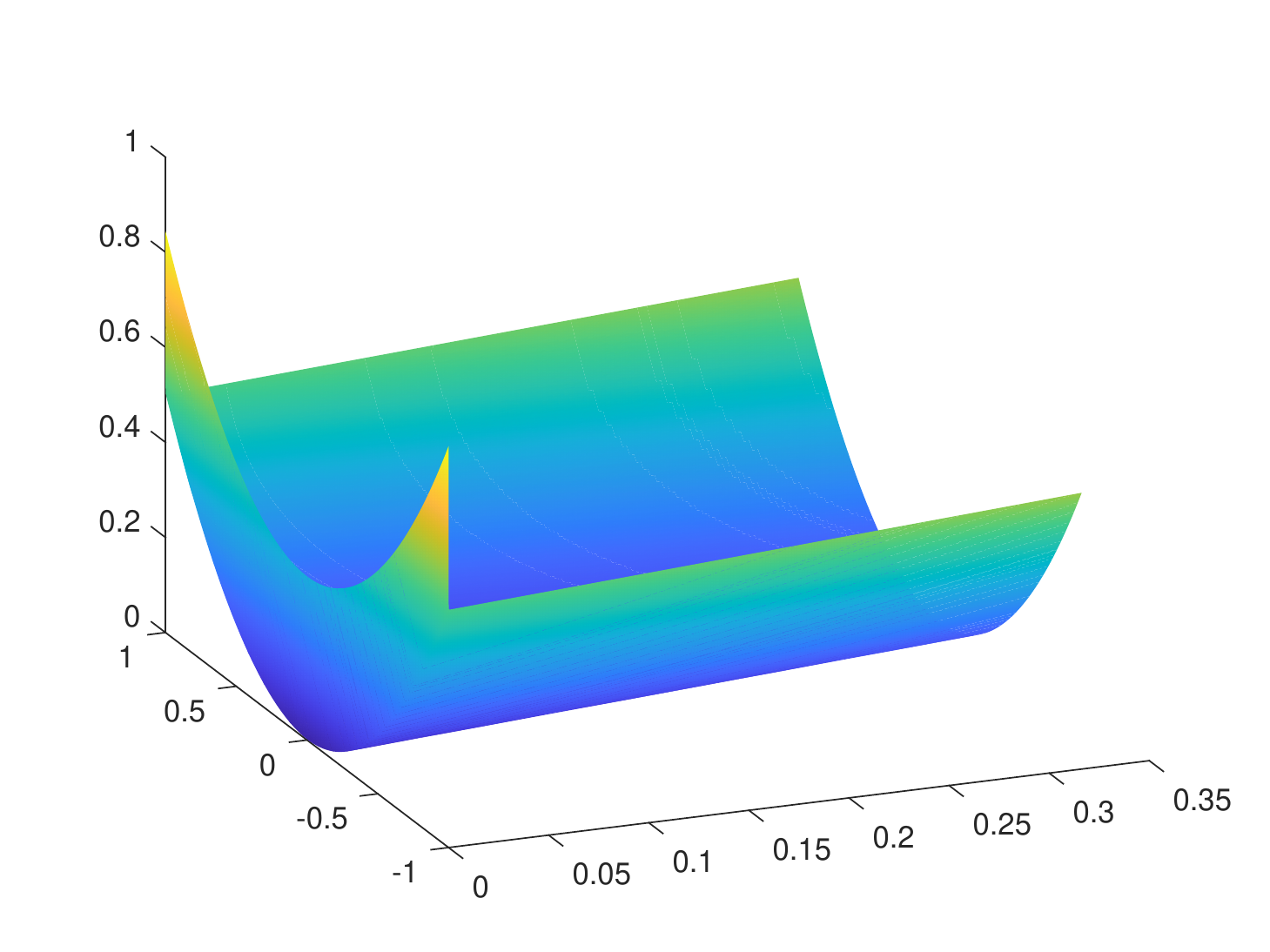}
			\put(60,1){$\phi$}
			\put(14,10){$\Phi$}
			\put(0,40){$\mathcal{H}$}
	\end{overpic}
	\caption{The Hamiltonian surface \cref{eq01:hamiltonian}.
 The sharp increase in $\mathcal{H}$ near $\phi=0$ is due to the repulsive singularity.}\label{fig01:hamiltonian_3D}
	\end{minipage}\quad\quad%
	\begin{minipage}{0.45\textwidth}
	\centering

	\begin{overpic}[percent,width=\textwidth]{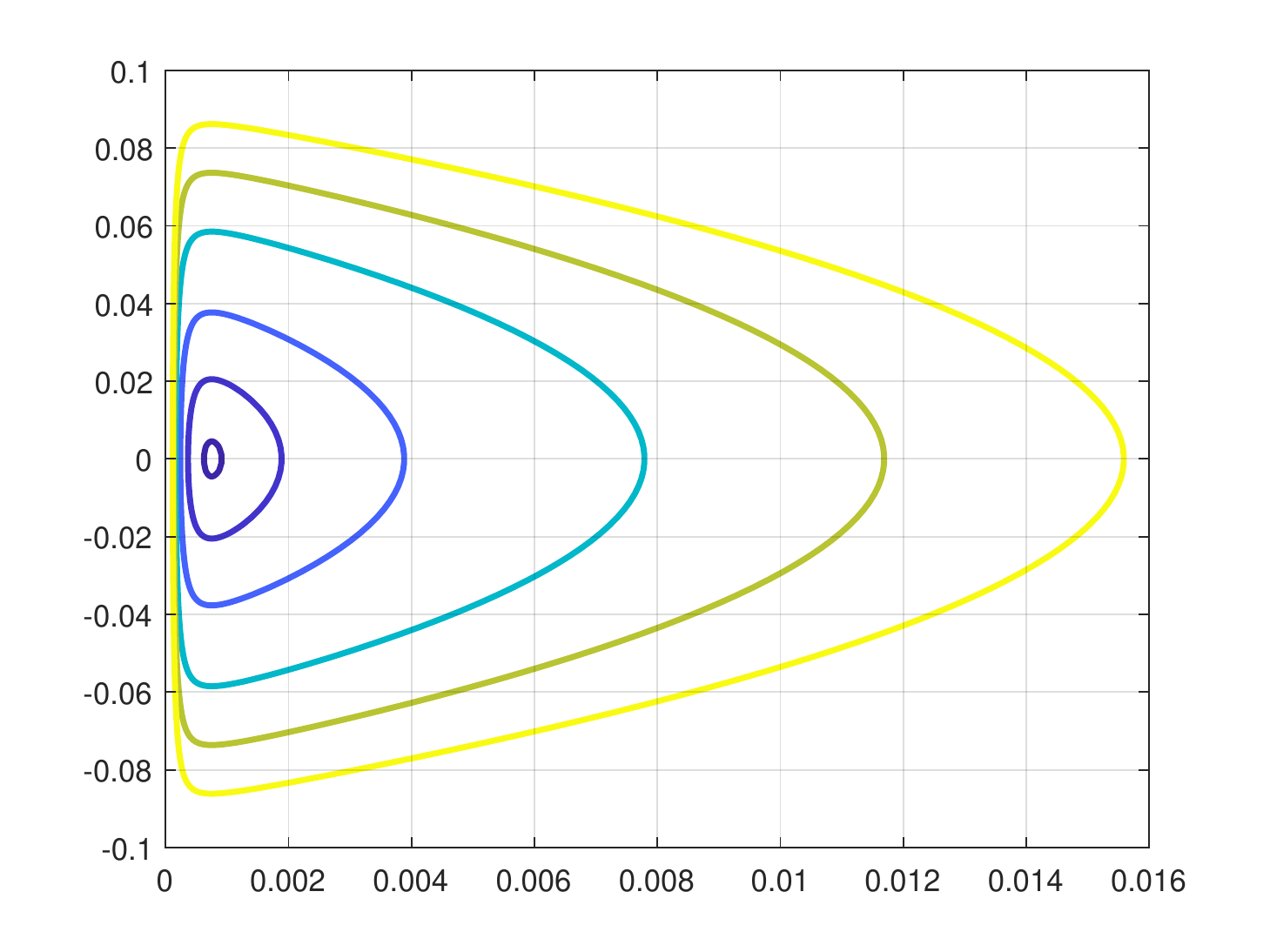}
			\put(50,0){$\phi$}
			\put(0,38){$\Phi$}
	\end{overpic}
	\caption{Contour plot of the Hamiltonian \cref{eq01:hamiltonian}, for $E = \num{7.5e-4},\num{1e-3}, \num{2e-3},\num{3e-3},\num{4e-3},\num{5e-3}$. Solutions of \cref{eq01:Phi_reduced} follow the contours in a clockwise direction.}\label{fig01:hamiltonian_contours}
	\end{minipage}
\end{figure}

\Cref{eq01:Phi_reduced} can be written as a planar Hamiltonian system with generalised coordinate $\phi$ and generalised momentum $\Phi$. Let $\mathcal{H}(\Phi,\phi)$ be the Hamiltonian
\begin{align}
	\mathcal{H}(\Phi,\phi) = & \frac{\Phi^2}{2} + \phi - \epsilon\Big(\ln^2(\phi)(\frac{a_{ll}}{2} + a_{ll}\phi) + (a_{l0} - 2a_{ll} - \frac{a_{l2}}{\phi})\ln(\phi) \\ \nonumber
 & + (a_0 -2a_{ll} - a_{l0})\phi - \frac{a_2 - a_{l2}}{\phi} - \frac{a_3}{2\phi^2}\Big). \label{eq01:hamiltonian1}
\end{align}
Then we recover \cref{eq01:Phi_reduced} by taking
\begin{align}
 	\diff{\Phi} =& -\frac{\partial \mathcal{H}}{\partial \phi} = \epsilon\left(\frac{a_3}{\phi^3} + \frac{a_{l2}\log\phi}{\phi^2} + \frac{a_{2}}{\phi^2} + \frac{a_{l1}\log\phi}{\phi} + \frac{a_{1}}{\phi} + a_{ll}(\log\phi)^2 + a_{l0}\log\phi + a_0\right)  - 1\label{eqA:ham1}\\
 	\diff{\phi} =& \frac{\partial \mathcal{H}}{\partial \Phi} = \Phi.  \label{eqA:ham2}
 \end{align}
The Hamiltonian surface $\mathcal{H}(\Phi,\phi)$ is shown in \cref{fig01:hamiltonian_3D}. 
Closed contours of the Hamiltonian indicate periodic orbits, seen in \cref{fig01:hamiltonian_contours}.
 An equilibrium point exists in the centre of the phase portrait, corresponding to the steady motion $S_\text{steady}$ \cref{eq01:equilibrium}. 
The symmetries in $\mathcal{H}(\Phi,\phi)$, clearly visible in \cref{fig01:hamiltonian_contours}, mean we can restrict our study to the upper half plane $\Phi>0$.

The period of one oscillation is given by 
\begin{align}
	T = 2\int^{\phi_1}_{\phi_0} \frac{1}{\Phi}\mathrm{d}\phi,
\end{align}
where $\phi_0$ and  $\phi_1$ are the intersections of the contour with the line $\Phi=0$. These intersections can be found from solutions of the Hamiltonian at $\Phi =0$ with a particular energy $E$
\begin{align}
	\mathcal{H}(0,\phi) = E
 \label{eq01:hamiltonian}
\end{align}
From \cref{fig01:hamiltonian_contours,fig01:hamiltonian_3D} we expect either two positive solutions for $\phi$, a single positive solution corresponding to the equilibrum or no solutions, depending on the value of $E$. 
In the case of $\epsilon=0$, 
\begin{align}
	T|_{\epsilon=0} = \sqrt{2}\int^{\phi_1}_{\phi_0}\frac{1}{\sqrt{h-\phi}}\mathrm{d}\phi.
\end{align}
If $\epsilon=0$, the can falls flat and $\phi_0=0$. If the can is released from rest at $\phi = \phi_1$, then $\phi_1$ is the highest angle attained by the can and $\mathcal{H} = \phi_1$. The period for such a motion is 
\begin{align}
	T|_{\epsilon=0} = 2\sqrt{2\phi_1}.
\end{align}
which agrees with the calculation for $\ddphi=-1$ given by  \cref{eq01:Phi_reduced} with $\epsilon=0$. In unscaled time 
\begin{align}
	T|_{\epsilon =0} = 2\sqrt{\frac{2\phi_1(A + mH^2 + mR^2)}{mgR}}.
\end{align}
As expected, cans with larger $H$ take a longer time to fall. While it is simpler to compute this leading order estimate of the oscillation period using \cref{eq01:Phi_reduced}, the Hamiltonian method of this section suggests a way to compute the period with non-zero $\epsilon$. 

\section{Coefficients in \texorpdfstring{\cref{eq01:order1_outer}}{(5.3)}}\label{app:b_coefficients}
\begin{align*}
    b_3 &= \frac{8a_3}{I^2}\\
    b_{l2} &= \frac{4a_{l2}}{I}\\
    b_2 &= \frac{4a_2}{I} + \frac{4a_{l2}\log(I/2)}{I}\\
    b_{l1} &= 2a_{l1}\\
    b_1 &= 2a_1 + 2a_{l1}\log(I/2)\\
    b_{ll} &= a_{ll}I\\
    b_{l0} &= a_{l0}I + 2a_{ll}I\log(I/2)\\
    b_0 &= a_0I + a_{l0}I\log(I/2) + a_{ll}I\log^2(I/2)
\end{align*}

\section{Integrals in \cref{eq01:integrals}} \label{secA:integrals}
The polynomial integrals
\begin{align}
	J_3 = \int\int \frac{1}{(1-T^2)^3}\mathrm{d}T\mathrm{d}T &= \frac{T^2}{8T^2-8} + \frac{3T}{8}\arctanh{T}\\
	J_2 = \int\int \frac{1}{(1-T^2)^2}\mathrm{d}T\mathrm{d}T &= \frac{T}{2}\arctanh(T) \\
	J_1 = \int\int \frac{1}{(1-T^2)}\mathrm{d}T\mathrm{d}T &= T\arctanh(T) +\frac{1}{2}\ln(1-T^2)
\end{align}
and the logarithmic integrals, computed in \textsc{Mathematica},
\begin{align}
J_{l2} = & \int \int \frac{\log \left(1-T^2\right)}{\left(1-T^2\right)^2} \, \mathrm{d}T \, \mathrm{d}T\\ \nonumber
=& \frac{1}{8} \Bigg[2 (T+1) \text{Li}_2\left(\frac{1-T}{2}\right)-2 (T-1) \text{Li}_2\left(\frac{T+1}{2}\right)-T \log ^2(1-T)+T \log ^2(T+1)\\
&-4 \left(\log (8) \log(1-T)+\log (T+1)+\log (8-8T)-3+2 \log ^2(2)\right)+4 \log (4) \log (1-T)\nonumber\\
&+2 \log (T+1) \log (1-T)+(T (\log (16)-4)-2 \log (4)) \tanh ^{-1}(T)\Bigg]\nonumber\\
J_{l1} =&\int \int \frac{\log \left(1-T^2\right)}{1-T^2} \, \mathrm{d}T\, \mathrm{d}T \\ \nonumber
=& \frac{1}{4} \Bigg[2 (T-1) \text{Li}_2\left(\frac{1-T}{2}\right)-2 (T+1) \text{Li}_2\left(\frac{T+1}{2}\right)-2 (T-1) \log ^2(1-T)+\log ^2(1-T)\\
&+(T+1) \log ^2(T+1)+(2 (T-1) \log (1-T)+\log (4)) \log (1-T)\nonumber\\
&-T \log (1-T) \log (4 (1-T))+(T+1) \log (4) \log (T+1)+8-2 \log (4)\Bigg]\nonumber\\
J_{l0}= &\int \int \log \left(1-T^2\right) \, \mathrm{d}T \, \mathrm{d}T \\ \nonumber
=& \frac{1}{2} \left(-3 T^2+\left(T^2+1\right) \log \left(1-T^2\right)+4 T \tanh ^{-1}(T)\right)\\
J_{ll}= &\int \int \log ^2\left(1-T^2\right) \, \mathrm{d}T \, \mathrm{d}T \\ \nonumber
=& \frac{1}{2} \Bigg[4 (T-1) \text{Li}_2\left(\frac{1-T}{2}\right)-4 (T+1) \text{Li}_2\left(\frac{T+1}{2}\right)+14 T^2\\
&+\log (1-T) \left(-6 T^2+2 \left(T^2-1\right) \log (T+1)+4 (T-1) \log (1-T)-2+\log (16)\right)\nonumber\\
&+(T-3) (T-1) \log ^2(1-T)-2 (T-1) \log ^2(1-T)+(T+1)^2 \log ^2(T+1)\nonumber\\
&-4 T (\log (2)-2) \log (1-T)-2 (T+1) (3 T+1-\log (4)) \log (T+1)+16-8 \log (2)\Bigg]\nonumber
\end{align}
where $\text{Li}_2$ is the dilogarithm.

\typeout{}
\bibliographystyle{plain}

\bibliography{library}


\end{document}